\documentclass[11pt, letter]{amsart}

\usepackage{amssymb,amscd, epsfig, mathrsfs, xypic, amsmath,amsthm,setspace}

\usepackage{ascmac}
\usepackage{fancybox}
\usepackage{color}
 
\xyoption{all}

\numberwithin{equation}{section} 
\theoremstyle{plain}
\newtheorem{theorem}{Theorem}[section]
\newtheorem{proposition}[theorem]{Proposition}         
\newtheorem{corollary}[theorem]{Corollary} 
\newtheorem{lemma}[theorem]{Lemma} 

\newtheorem{remark}[theorem]{Remark}

\renewcommand{\proof}{\no \textit{Proof}$\ \\ \no$}

\newcommand{\no}{\noindent}

\newcommand{\C}{\mathbb C}   
\newcommand{\R}{\mathbb R}

\newcommand{\Z}{\mathbb Z}

\newcommand{\fcw}{\omega}
\newcommand{\dPhi}{\Phi^*}
\newcommand{\ddPhi}{\Pi^*}
\newcommand{\RL}{M}
\newcommand{\CL}{N}
\newcommand{\SimR}{\Pi}
\newcommand{\rnk}{n}
\newcommand{\coi}{i}

\newcommand{\parafir}{c}
\newcommand{\parasec}{a}
\newcommand{\parathr}{b}
\newcommand{\parafor}{y}
\newcommand{\A}{S}
\newcommand{\intnum}[1]{(\mu_X,{#1})}
\newcommand{\invdiv}[1]{D_{#1}}
\newcommand{\Asetfir}{A}
\newcommand{\Asetsec}{B}
\newcommand{\Asetthr}{C}
\newcommand{\Asetrun}{B'}
\newcommand{\Bsetfir}{A}
\newcommand{\Bsetsec}{B}
\newcommand{\Bsetthr}{C}
\newcommand{\Bsetrun}{B'}

\newcommand{\Csetsec}{B}

\newcommand{\Csetrun}{B'}
\newcommand{\Dsetfir}{A}
\newcommand{\Dsetsec}{B}
\newcommand{\Dsetthr}{C}
\newcommand{\Dsetrun}{B'}
\newcommand{\DA}{S^{\circ}}
\newcommand{\Euc}{E}
\newcommand{\WY}[3]{\lambda_{#1 #2}^{#3}}

\newcommand{\SC}[3]{c_{#1 #2}^{#3}}
\newcommand{\imat}{\mathcal{I}}
\newcommand{\imatcomp}[2]{\mathcal{I}^{#1}_{\ #2}}
\newcommand{\imatcompinv}[2]{(\mathcal{I}^{-1})^{#1}_{\ #2}}
\newcommand{\WG}{W_{G_2}}

\allowdisplaybreaks[1]

\begin{document}
\title[Young diagrams and intersection numbers]{Young diagrams and intersection numbers for \\ toric manifolds associated with Weyl chambers}
\author{Hiraku Abe}
\address{Osaka City University Advanced Mathematical Institute, Sumiyoshi-ku, Osaka 558-8585, Japan}
\email{hirakuabe@globe.ocn.ne.jp}

\date{\today}
\keywords{Young diagrams, intersection numbers, toric varieties, structure constants} \subjclass[2010]{
14M25, 
17B22, 
13F55. 
}

\maketitle

\begin{abstract}
We study intersection numbers of invariant divisors in the toric manifold associated with the fan determined by the collection of Weyl chambers for each root system of classical type and of exceptional type $G_2$. 
We give a combinatorial formula for intersection numbers of certain subvarieties which are naturally indexed by elements of the Weyl group. These numbers describe the ring structure of the cohomology of the toric manifold. 
\end{abstract}

\section{Introduction}\label{introduction}
Let $\Phi$ be a root system in the $\rnk$-dimensional Euclidean space $\Euc$ with its inner product. We denote by $\Delta(\Phi)$ the fan determined by the collection of Weyl chambers in $E^*$, and consider the toric manifold $X$ associated with $\Delta(\Phi)$.
This toric manifold arises as the closure of a general orbit in the flag variety with respect to the standard torus action which makes $X$ a regular semisimple Hessenberg variety (\cite{De Mari-Procesi-Shayman}). 
The Weyl group $W$ naturally acts on the Weyl chambers and hence also on $X$.
The representation of $W$ on the cohomology $H^*(X;\C)$ has been extensively studied by Procesi \cite{Procesi}, Dolgachev-Lunts \cite{Dolgachev-Lunts}, and Stembridge \cite{Stembridge}. 
For the classical root system of type $A_{\rnk}$, Losev-Manin \cite{Losev-Manin} described $X$ as the moduli space of stable $(\rnk+1)$-pointed chains of projective lines (cf. Batyrev-Blume \cite{Batyrev-Blume}).

Let $\Pi=\{\alpha_1,\cdots,\alpha_{\rnk}\}\subset \Phi$ be a set of simple roots, then we have a torus invariant non-singular subvariety $X_u$ in $X$ of codimension $|u(\Pi)\cap\Phi^-|$ such that the associated cohomology classes $\{[X_u]\}_{u\in W}$ form a module basis of the integral singular cohomology $H^*(X)$.
The cohomology class $[X_u]$ is written as a monomial of torus invariant divisors $\invdiv{u\omega_i}$ of $X$ for all coweights $u\omega_i$ satisfying $u\alpha_i\in\Phi^-$ where $\{\omega_1,\cdots,\omega_{\rnk}\}$ is the set of fundamental coweights (see Section \ref{prelim} and (\ref{def of Xu}) for details).
In this paper, we study the case for the root systems of classical type and of exceptional type $G_2$, and we give a combinatorial formula of the intersection numbers 
\begin{align*}
\intnum{[w_0X_{w_0 w}][X_u][X_v]}
\end{align*} 
of three subvarieties $X_u$, $X_v$ and $w_0X_{w_0 w}$ for $u,v,w\in W$
where $\mu_X$ is the fundamental homology class of $X$ and
$w_0$ is the longest element.
As an application, 
we will obtain a recursive formula for the structure constants $c_{u,v}^w$ in the expansion of the product
\begin{align*}
[X_u][X_v] = \sum_{w\in W} c_{u,v}^w [X_w]
\quad \text{where} \quad c_{u,v}^w \in\Z
\end{align*} 
as discussed in Section 4.

\vspace{10pt}
Let us state our formula for $\intnum{[w_0X_{w_0 w}][X_u][X_v]}$ in the case of the classical root system of type $A_{\rnk}$ (the results for the classical root systems of type $B_{\rnk}, C_{\rnk}$, and $D_{\rnk}$ are stated in Section 5). 
In this case, the Weyl group $W$ is the $(\rnk+1)$-th permutation group $\mathfrak{S}_{\rnk+1}$.
For each $u\in \mathfrak{S}_{\rnk+1}$, we let 
\begin{align}\label{def of D(u)}
&D(u):=\{ \{u(1),u(2),\cdots,u(i)\} \mid u(i)>u(i+1) \},  \\ \label{def of A(u)}
&A(u):=\{ \{u(1),u(2),\cdots,u(i)\} \mid u(i)<u(i+1) \}
\end{align}
where each $\{u(1),u(2),\cdots,u(i)\}$ is a subset of $[\rnk+1]$.
We define a Young diagram $\WY{u}{v}{w}$ as follows. Assume $d(u)+d(v)=d(w)$, and the collection $D(u)\coprod D(v) \coprod A(w)$ forms a nested chain of subsets of $[\rnk+1]$. In this case, $\WY{u}{v}{w}$ is defined to be the Young diagram consisting of the cardinalities of those chains ordered as a weakly decreasing sequence. Otherwise, $\WY{u}{v}{w}=\emptyset$.
For example, suppose $n=4$ and, let $u=12354$, $v=31254$, and $w=35421$. Then, we have that $D(u)=\{1235\}$, $D(v)=\{3, 3125\}=\{3, 1235\}$, and $A(w)=\{3\}$ where $1235$ denotes the set $\{1,2,3,5\}$ and we use the same notation for others. These sets forms a nested chain of subsets $3\subset 3\subset 1235\subset 1235$, and hence we obtain $\lambda_{12354, 31254}^{35421}=(4,4,1,1)$.

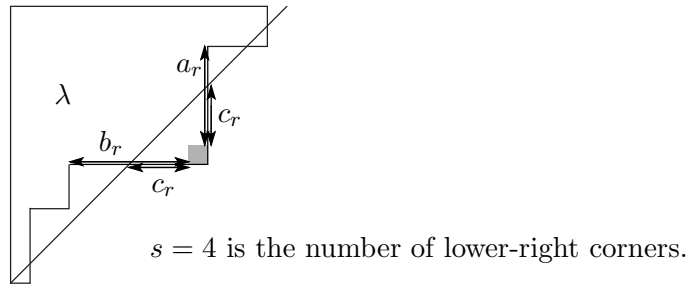
\begin{figure}[h]
\centering
\unitlength 0.1in
\begin{picture}( 30.0300, 16.0900)( 16.0000,-24.8900)
%
\special{pn 8}%
\special{pa 1704 2418}%
\special{pa 1806 2418}%
\special{fp}%
%
\special{pn 8}%
\special{pa 1806 2418}%
\special{pa 1806 2026}%
\special{fp}%
%
\special{pn 8}%
\special{pa 1806 2028}%
\special{pa 2008 2028}%
\special{fp}%
%
\special{pn 8}%
\special{pa 2010 2028}%
\special{pa 2010 1798}%
\special{fp}%
%
\special{pn 8}%
\special{pa 2736 1178}%
\special{pa 3048 1178}%
\special{fp}%
%
\special{pn 8}%
\special{pa 3048 1178}%
\special{pa 3048 968}%
\special{fp}%
%
\special{pn 8}%
\special{pa 3048 968}%
\special{pa 1704 968}%
\special{fp}%
%
\special{pn 8}%
\special{pa 1704 968}%
\special{pa 1704 2418}%
\special{fp}%
%
\special{pn 8}%
\special{pa 1704 2418}%
\special{pa 3150 968}%
\special{fp}%
%
\special{pn 8}%
\special{pa 2754 1488}%
\special{pa 2754 1694}%
\special{fp}%
\special{sh 1}%
\special{pa 2754 1694}%
\special{pa 2774 1626}%
\special{pa 2754 1640}%
\special{pa 2734 1626}%
\special{pa 2754 1694}%
\special{fp}%
\special{pa 2754 1694}%
\special{pa 2754 1384}%
\special{fp}%
\special{sh 1}%
\special{pa 2754 1384}%
\special{pa 2734 1452}%
\special{pa 2754 1438}%
\special{pa 2774 1452}%
\special{pa 2754 1384}%
\special{fp}%
%
\special{pn 8}%
\special{sh 0.300}%
\special{pa 2636 1694}%
\special{pa 2736 1694}%
\special{pa 2736 1796}%
\special{pa 2636 1796}%
\special{pa 2636 1694}%
\special{ip}%
%
\special{pn 8}%
\special{pa 2736 1796}%
\special{pa 2736 1178}%
\special{fp}%
%
\special{pn 8}%
\special{pa 2720 1694}%
\special{pa 2720 1178}%
\special{fp}%
\special{sh 1}%
\special{pa 2720 1178}%
\special{pa 2700 1246}%
\special{pa 2720 1232}%
\special{pa 2740 1246}%
\special{pa 2720 1178}%
\special{fp}%
\special{pa 2720 1178}%
\special{pa 2720 1694}%
\special{fp}%
\special{sh 1}%
\special{pa 2720 1694}%
\special{pa 2740 1626}%
\special{pa 2720 1640}%
\special{pa 2700 1626}%
\special{pa 2720 1694}%
\special{fp}%
%
\special{pn 8}%
\special{pa 2016 1782}%
\special{pa 2636 1782}%
\special{fp}%
\special{sh 1}%
\special{pa 2636 1782}%
\special{pa 2568 1762}%
\special{pa 2582 1782}%
\special{pa 2568 1802}%
\special{pa 2636 1782}%
\special{fp}%
\special{pa 2636 1782}%
\special{pa 2016 1782}%
\special{fp}%
\special{sh 1}%
\special{pa 2016 1782}%
\special{pa 2082 1802}%
\special{pa 2068 1782}%
\special{pa 2082 1762}%
\special{pa 2016 1782}%
\special{fp}%
%
\special{pn 8}%
\special{pa 2324 1812}%
\special{pa 2636 1812}%
\special{fp}%
\special{sh 1}%
\special{pa 2636 1812}%
\special{pa 2568 1792}%
\special{pa 2582 1812}%
\special{pa 2568 1832}%
\special{pa 2636 1812}%
\special{fp}%
\special{pa 2636 1812}%
\special{pa 2324 1812}%
\special{fp}%
\special{sh 1}%
\special{pa 2324 1812}%
\special{pa 2392 1832}%
\special{pa 2378 1812}%
\special{pa 2392 1792}%
\special{pa 2324 1812}%
\special{fp}%
%
\special{pn 8}%
\special{pa 2016 1796}%
\special{pa 2736 1796}%
\special{fp}%
\put(27.8800,-15.9600){\makebox(0,0)[lb]{$\parafir_{r}$}}%
\put(25.7000,-13.3300){\makebox(0,0)[lb]{$\parasec_{r}$}}%
\put(21.6700,-17.4600){\makebox(0,0)[lb]{$\parathr_{r}$}}%
\put(24.4200,-19.6000){\makebox(0,0)[lb]{$\parafir_{r}$}}%
\put(24.3400,-23.1200){\makebox(0,0)[lb]{$s=4$ is the number of lower-right corners.}}%
\put(19.3400,-14.8000){\makebox(0,0)[lb]{$\lambda$}}%
%
\special{pn 8}%
\special{pa 1600 880}%
\special{pa 4604 880}%
\special{pa 4604 2490}%
\special{pa 1600 2490}%
\special{pa 1600 880}%
\special{ip}%
\end{picture}%
\vspace{-5pt}
\caption{Three numbers $\parasec_r$, $\parathr_r$ and $\parafir_r$}
\label{Intro picture}
\end{figure}

For a Young diagram $\lambda=(\lambda_1\geq\cdots\geq\lambda_{\rnk})$ with $\rnk$ rows (i.e. $\lambda_{\rnk}>0$) fitting into the $\rnk\times\rnk$ square, we define $I(\lambda)\in\Z$ to be the following integer.
Let $s$ be the number of lower-right corners of $\lambda$, i.e.,  $s=|\{ i\in[\rnk] \mid \lambda_{i}>\lambda_{i+1} \}|$ where $\lambda_{\rnk+1}:=0$.
Write
\begin{align*}
\{ i\in[\rnk] \mid \lambda_{i}>\lambda_{i+1} \} 
=\{ \coi_1,\cdots,\coi_{s} \}.
\end{align*}
We impose the condition $\coi_1<\coi_2<\cdots<\coi_{s}$ to determine them uniquely. 
Observe that $\coi_{s}=\rnk$.
For example, if $\rnk=4$ and $\lambda=(4,2,2,1)$, then $s=3$ and $\{\coi_1,\coi_2, \coi_3\}=\{1,3,4\}$.
For $r=1,\cdots,s$, define
\begin{align*}
\parasec_{r}:=\coi_{r}-\coi_{r-1}-1, \quad
\parathr_{r}:=\lambda_{i_{r}}-\lambda_{i_{r+1}}-1, \quad
\parafir_{r}:=\lambda_{i_{r}}+\coi_{r}-\rnk-1
\end{align*}
where we write $\coi_0=0$, and let
\begin{align*}
\parafor_r:=
\binom{ \parasec_{r} }{ \parafir_{r} }
\binom{ \parathr_{r} }{ \parafir_{r} }
\quad \text{for $r=1,\cdots,s$}.
\end{align*}
We use the convention $\binom{x}{y}=0$ unless $0\leq y\leq x$.
By shading each lower-right corner of $\lambda$, the pictorial meanings of $\parasec_r$, $\parathr_r$, and $\parafir_r$ become clear (as shown in Figure \ref{Intro picture}).
Namely, $\parasec_r$ is the number of boxes between the north side of the shaded box and the the upper-left corner placed above,
$\parathr_r$ is the similar number for the horizontal segment of the corner, and $\parafir_r$ is the number of boxes between the north side of the shaded box and the crossing point of the vertical segment and the anti-diagonal line where we count negatively if that part of the vertical segment is above the anti-diagonal. 
\vspace{-5pt}

\noindent
Now, let 
\begin{align*}
I(\lambda) := (-1)^{\rnk+s}\parafor_{1}\cdots \parafor_{s},
\end{align*}
and put $I(\emptyset)=0$.
The following is our main statement for type $A_{\rnk}$.
\begin{theorem}\label{intro main thm}
$\displaystyle{\intnum{[w_0X_{w_0w}][X_u][X_v]}=I(\WY{u}{v}{w})}$
for any $u,v,w\in\mathfrak{S}_{\rnk+1}$ 
where $\mu_X$ is the fundamental homology class of $X$.
\end{theorem}

We will prove Theorem \ref{intro main thm} by computing general intersection numbers of invariant divisors of $X$ in Section 3 and 4.
Section 5 is devoted to the classical root systems of type $B_{\rnk}$, $C_{\rnk}$ and $D_{\rnk}$.

\subsection*{Acknowledgements}
The author would like to thank Tatsuya Horiguchi, Hiroaki Ishida, Ivan Limonchenko and Tomoo Matsumura for valuable comments. He also thank Miho Hatanaka for reading of the first version of this paper carefully.

\vspace{20pt}
\section{Preliminaries}\label{prelim}
Let $\Phi$ be a root system in the $\rnk$-dimensional Euclidean space $\Euc$ with its inner product. Let $\RL\subset \Euc$ be the root lattice of $\Phi$ and $\CL\subset \Euc^*$ be the coweight lattice of $\Phi$.  Then $\RL$ is the dual lattice of $\CL$ with respect to the natural pairing.

We choose a set of simple roots $\SimR=\{\alpha_1,\cdots,\alpha_{\rnk}\}\subset\Phi\subset \Euc$, and let $\ddPhi:=\{\fcw_1,\cdots,\fcw_{\rnk}\}\subset E^*$ be the dual basis of $\SimR$ defined by $\langle \fcw_i, \alpha_j \rangle=\delta_{ij}$, i.e., $\fcw_1,\cdots,\fcw_{\rnk}$ are the \textit{fundamental coweights}.
For each $u\in W$, denote 
\begin{align*}
\sigma_{u}
:=\text{cone}(u\fcw_1,\cdots,u\fcw_{\rnk}) 
= \{ \textstyle{\sum_{i=1}^{\rnk}} \lambda_iu\fcw_i \mid \lambda_i\geq0 \}.
\end{align*} 
These cones $\{\sigma_{u}\}_{u\in W}$ form a non-singular complete fan $\Delta(\Phi)$ in $\Euc^*$ by including all their faces. 
The set of minimal generators of these cones are the set of \textit{coweights}:
\begin{align*}
 \dPhi=\bigcup_{v\in W} \{v\fcw_1,\cdots,v\fcw_{\rnk}\}.
\end{align*} 
For each element $u\in W$, 
the maximal cones containing a minimal generator $u\fcw_i$ are $\sigma_v$ for $v\in W$ such that $u\fcw_i=v\fcw_j$ for some $j$.
There is a cone of $\Delta(\Phi)$ generated by minimal generators $x_1,\cdots,x_k\in\dPhi$ if and only if there exists $u\in W$ such that each $x_i$ can be written as $u\fcw_{j}$ for some $j$.

We consider the toric manifold $X=X(\Phi)$ associated with the fan $\Delta(\Phi)$. 
For root systems $\Phi$ and $\Phi'$, it is easily verified that $X(\Phi)\cong X(\Phi')$ as toric varieties (in the sense of \cite{Cox-Little-Schenck} Sec. 3.3.) if and only if $\Phi\cong \Phi'$ as root systems (in the sense that their Cartan matrices are the same up to permuting the indexes).
We refer to \cite{Batyrev-Blume} and \cite{Klyachko} for general properties of $X$. 
Let $\invdiv{x}\subset X$ be the invariant divisor corresponding to the ray generated by $x\in\dPhi$. 
The Poincar$\acute{\text{e}}$ dual $\tau_{x}:=[\invdiv{x}]$ gives us a cohomology class of degree $2$ in the integral singular cohomology $H^*(X)$.
%
The cohomology ring $H^*(X)$ is isomorphic to the face ring of the underlying simplicial complex of the fan $\Delta(\Phi)$ modulo some linear relations (\cite{Fulton1}). More precisely, we have 
\begin{align*}
H^*(X) = \Z[\tau_x \mid x\in\dPhi]/I
\end{align*} 
where the ideal $I$ is generated by $\tau_{x_1}\cdots\tau_{x_k}$ for which $x_1,\cdots,x_k$ do not generate a face of  $\sigma_u$ for some $u\in W$ and $\sum_{x\in\dPhi} \langle x,\alpha \rangle \tau_{x}$ for any root $\alpha$.
Namely, we have the following equalities in $H^*(X)$ :
\begin{align}\label{prelim linear relation}
\sum_{x\in\dPhi} \langle x,\alpha \rangle \tau_{x}=0
\quad \text{for any root $\alpha$}.
\end{align} 

The above observation about rays of $\Delta(\Phi)$ implies that 
\begin{lemma}\label{prelim rays generating a cone}
We have $\tau_{x_1}\cdots\tau_{x_k} = 0$ unless there exists $u\in W$ such that each $x_i$ can be written as $u\fcw_{j}$ for some $j$.
\end{lemma}

Let $\mu_X$ be the fundamental homology class of $X$.
For subvarieties $Z_1,\cdots,Z_k\subset X$, we call $\intnum{[Z_1]\cdots[Z_k]}$ the \textit{intersection number} of $Z_1,\cdots,Z_k$ where $[Z_i]$ denotes the Poincar$\acute{\text{e}}$ dual of $Z_i$.
Note that the Weyl group $W$ acts on the fan $\Delta(\Phi)$, and hence acts on the toric manifold $X$.
We have $uX_{x}=X_{ux}$ for any $u\in W$ and $x\in\dPhi$ which means that $(u^{-1})^*\tau_{x} = \tau_{ux}$.
The next lemma says that intersection numbers for divisors $\invdiv{u\fcw_i}$ are invariant under the Weyl group action. 
\begin{lemma}\label{prelim Weyl inv}
Let $x_1,\cdots,x_{\rnk}\in\dPhi$. Then for any $u\in W$, we have
\begin{align*}
\intnum{\tau_{ux_1}\cdots\tau_{ux_{\rnk}}} = \intnum{\tau_{x_1}\cdots\tau_{x_{\rnk}}}.
\end{align*} 
\end{lemma}
\proof
Observe that $\tau_{ux_1}\cdots\tau_{ux_{\rnk}}=(u^{-1})^*(\tau_{x_1}\cdots\tau_{x_{\rnk}})$. Both of $\tau_{ux_1}\cdots\tau_{ux_{\rnk}}$ and $\tau_{x_1}\cdots\tau_{x_{\rnk}}$ can be written as the cohomology class $[p]$ of a point $p$ in $X$ multiplied by some integer, and these integers are the corresponding intersection numbers.
Since $u$ preserves the orientation of $X$, we have $(u^{-1})^*([p])=[u\cdot p]=[p]$ which proves the claim.
\qed

\vspace{10pt}
For any $u\in W$, the product $\tau_{u\fcw_{1}} \cdots \tau_{u\fcw_{\rnk}}$ is exactly the Poincar$\acute{\text{e}}$ dual of a point in $X$ since the invariant divisors $X_{u\fcw_1},\cdots,X_{u\fcw_{\rnk}}$ intersect transversally which means that
\begin{align}\label{general type I=1}
\intnum{\tau_{u\fcw_{1}} \cdots \tau_{u\fcw_{\rnk}}}=1.
\end{align}
For the root systems of classical type, we will compute the intersection number $\intnum{\tau_{x_1}\cdots\tau_{x_{\rnk}}}$ for arbitrary $x_1,\cdots,x_{\rnk}\in\dPhi$.
By Lemma \ref{prelim rays generating a cone}, we can assume that this number is of the form $\intnum{(\tau_{u\fcw_{i_1}})^{m_1} \cdots (\tau_{u\fcw_{i_s}})^{m_s}}$ for some $1\leq i_1< \cdots< i_s\leq\rnk$ and $1\leq m_{k}\leq\rnk \ (k=1,\cdots,s)$ satisfying $m_1+\cdots+m_s=n$ without loss of generality.
We call the number $m_{k}$ the \textit{multiplicity} of $\tau_{u\omega_{i_k}}$.
We compute this number by applying the linear relations (\ref{prelim linear relation}) to reduce the multiplicities $m_1,\cdots,m_s$. 
Although Lemma \ref{prelim Weyl inv} shows that this number is equal to $\intnum{(\tau_{\fcw_{i_1}})^{m_1} \cdots (\tau_{\fcw_{i_s}})^{m_s}}$, we will need Lemma \ref{prelim Weyl inv} again after applying the relations (\ref{prelim linear relation}).

In the next section, we will consider the classical root system of type $A_{\rnk}$, and compute the intersection numbers.

\vspace{10pt}
\section{Intersection numbers for Type $A_{\rnk}$}\label{section type A}
In this section, we compute the intersection numbers for the toric manifold $X$ of type $A_{\rnk}$. 
Let $E=\{x\in\R^{\rnk+1} \mid x_1+\cdots+x_{\rnk+1}=0\}$.
The roots are  $t_i-t_j \in E \ (1\leq i,j\leq \rnk+1)$ where $t_i\in\R^{\rnk+1}$ is the $i$-th standard vector.
We choose $\SimR=\{t_i-t_{i+1}\mid 1\leq i\leq \rnk\}$
as the set of simple roots, and write $\alpha_i=t_i-t_{i+1}$ for each $i$.
The Weyl group $W=\mathfrak{S}_{\rnk+1}$ is the $(\rnk+1)$-th permutation group acting on $E$ by $u(t_i-t_j)=t_{u(i)}-t_{u(j)}$ for each $u\in W$.
The minimal generators $\fcw_1,\cdots,\fcw_{\rnk}\in E^*$ of the fundamental Weyl chamber $\sigma_{\text{id}}$ are
\[
\fcw_i=(e_1+\cdots+e_i)-\frac{i}{\rnk+1}(e_1+\cdots+e_{\rnk+1}) 
\quad \text{for} \quad
i=1,\cdots,\rnk
\]
where $\{e_i\}_i\subset (\R^{\rnk+1})^*$ is the dual basis of $\{t_i\}_i\subset \R^{\rnk+1}$.

Denoting by $2^{[\rnk+1]}$ the set of all subsets of $[\rnk+1]=\{1,\cdots,\rnk+1\}$, we have a well-defined map $\dPhi \rightarrow 2^{[\rnk+1]}$ by sending $u\fcw_i \mapsto \{u(1),\cdots,u(i)\}$.
It is easy to see that this is an injection, and this establishes an identification
\begin{align}\label{type A identification}
\dPhi \quad \longleftrightarrow \quad \text{the set of non-empty proper subsets of $[\rnk+1]$}.
\end{align} 
In particular, the well-definedness implies that if $u\fcw_i=v\fcw_j$ then $i=j$.
Now, for each $\emptyset\subsetneq S\subsetneq[\rnk+1]$, we define
$\tau_S:=\tau_{u\omega_i}$ where $u\omega_i\in\dPhi$ corresponds to $S$ by this identification. 
Then, for $\emptyset\subsetneq S_1, \cdots, S_q\subsetneq[\rnk+1] \ (1\leq q\leq\rnk)$, it follows by Lemma \ref{prelim rays generating a cone} that $\tau_{S_1}\cdots\tau_{S_q}=0$ unless these sets form a nested chain of subsets, i.e. $S_1\subset\cdots\subset S_q$ up to reordering.

With Lemma \ref{prelim Weyl inv}, it is easy to show the following invariance property of intersection numbers which implies that $\intnum{\tau_{S_1}\cdots\tau_{S_{\rnk}}}$ for $\emptyset\subsetneq S_1\subset \cdots \subset S_{\rnk} \subsetneq [\rnk+1]$ is determined by the set of integers $1\leq |S_1|\leq \cdots\leq |S_{\rnk}|\leq \rnk$.
\begin{lemma}\label{invariance for type A}
Let $\emptyset\subsetneq S_1\subset \cdots \subset S_{\rnk} \subsetneq [\rnk+1]$ and $\emptyset\subsetneq S'_1\subset \cdots \subset S'_{\rnk} \subsetneq [\rnk+1]$.
If $|S_i|=|S'_i|$ for all $i=1,\cdots,\rnk$, then $\intnum{\tau_{S_1}\cdots\tau_{S_{\rnk}}} = \intnum{\tau_{S'_1}\cdots\tau_{S'_{\rnk}}}.$
\end{lemma}

Motivated by this property, we compute intersection numbers in terms of Young diagrams consisting of the cardinalities of the sets corresponding to the given invariant divisors. 
The linear relations (\ref{prelim linear relation}) are translated to 
\begin{align}\label{SR relation for type A}
\sum_{\substack{\emptyset \subsetneq S\subsetneq [\rnk+1] \\ k\in S, l\notin S}} \tau_{S} 
-
\sum_{\substack{\emptyset \subsetneq S\subsetneq [\rnk+1] \\ k\notin S, l\in S}} \tau_{S} 
=0
\qquad  \text{for each $k, l\in[\rnk+1]$}.
\end{align} 
In the following, we write $\tau_{\emptyset}=\tau_{[\rnk+1]}=1$.
This equality together with the above observation about $\tau_{S_1}\cdots\tau_{S_q}$ being $0$ implies the next lemma.
\begin{lemma}\label{type A separation}
Let $\emptyset\subset \Asetfir\subsetneq \Asetsec \subsetneq \Asetthr \subset [\rnk+1]$.
For any $b\in \Asetsec\backslash \Asetfir$ and $c\in \Asetthr\backslash \Asetsec$, we have
\begin{align*}
\tau_{\Asetfir} {\tau_{\Asetsec}}^2 \tau_{\Asetthr} = - \sum_{\substack{\Asetfir\subsetneq \Asetrun \subsetneq \Asetthr \\ \Asetrun\neq \Asetsec, b\in \Asetrun, c\notin \Asetrun}} \tau_{\Asetfir} \tau_{\Asetsec}\tau_{\Asetrun} \tau_{\Asetthr}.
\end{align*}
\end{lemma}

\vspace{10pt}
For a Young diagram fitting into the $\rnk\times\rnk$ square, we write the \textit{dotted anti-diagonal line} shifted down half the length of a single box from the standard anti-diagonal. 

\begin{figure}[h]
\centering
\input{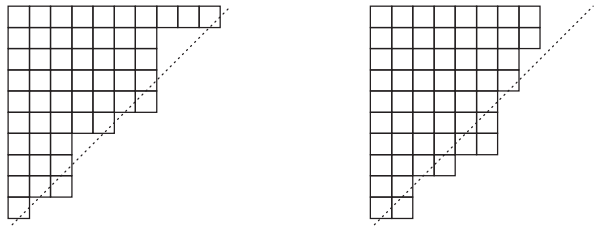}
\caption{Young diagrams and the dotted anti-diagonal line}
\label{zigzag picture}
\end{figure}

Let $\emptyset\subsetneq S_1\subset \cdots \subset S_{\rnk} \subsetneq [\rnk+1]$, and denote by $\lambda$ the Young diagram consisting of $\lambda_i=|S_{\rnk+1-i}|$ for each $i$.
We write $\lambda_{\rnk+1}=0$.
Let $s$ be the number of the lower-right corners of $\lambda$, that is, 
\begin{align*}
s:=|\{ i\in[\rnk] \mid \lambda_{i}>\lambda_{i+1} \}|.
\end{align*}
\begin{proposition}\label{dotted vanishing type A}
\emph{(Vanishing property)}
$\intnum{\tau_{S_1}\cdots\tau_{S_{\rnk}}} = 0$ unless each step of the zigzag line of the lower-right corners of $\lambda$ crosses the dotted anti-diagonal. 
\end{proposition}
\proof
We suppose that there is a step of the zigzag line of $\lambda$ which does not cross the dotted anti-diagonal, and show  $\intnum{\tau_{S_1}\cdots\tau_{S_{\rnk}}} = 0$ by induction on $k:=\rnk-s$.
Since there is no such case for $k=0$, we consider the case $k=1$. In this case,
there is a unique vertical segment of length 2 in the zigzag line of $\lambda$.
If there is a (unique) horizontal segment of length 2, then the vertical and horizontal segments are not adjacent because of our assumption.
By applying Lemma \ref{type A separation} for the square corresponding to this vertical segment, it follows that the intersection number is zero since there is no summand.

For the general case, take a vertical segment of length $\geq2$.
Let us say that this vertical segment contains $S_{i}$ and $S_{i+1}$ (i.e. $S_{i}=S_{i+1}$). 
We separate this square in $\tau_{S_1}\cdots\tau_{S_{\rnk}}$ by Lemma \ref{type A separation}. Let $\lambda'$ be the Young diagram corresponding to a summand in the right-hand-side. 
Then the zigzag line of $\lambda'$ has a step which does not cross the dotted anti-diagonal. 
In fact, if the vertical segment does not cross the dotted anti-diagonal, then this segment survives as a non-crossing segment of length at least 1, and if it does then we can find another vertical segment which does not, and this segment is preserved for each $\lambda'$ in the summands.
Now the induction hypothesis shows that each term will vanish after taking the intersection number, and we get $\intnum{\tau_{S_1}\cdots\tau_{S_{\rnk}}} = 0$.
\qed

\vspace{20pt}
Let $\lambda=(\lambda_1\geq\cdots\geq\lambda_{\rnk})$ be a Young diagram with $\rnk$ rows (i.e. $\lambda_{\rnk}>0$) fitting into the $\rnk\times\rnk$ square.
Let $I(\lambda)\in\Z$ be the one defined in Section 1.
We here recall the definition for the convenience of the reader.
Let $s$ be the number of lower-right corners of $\lambda$, i.e.,  $s=|\{ i\in[\rnk] \mid \lambda_{i}>\lambda_{i+1} \}|$ where $\lambda_{\rnk+1}:=0$.
Write
\begin{align*}
\{ i\in[\rnk] \mid \lambda_{i}>\lambda_{i+1} \} 
=\{ \coi_1,\cdots,\coi_{s} \}.
\end{align*}
We impose the condition $\coi_1<\coi_2<\cdots<\coi_{s}$ to determine them uniquely. 
Observe that $\coi_{s}=\rnk$.
For $r=1,\cdots,s$, define
\begin{align}\label{definition of a b c}
\parasec_{r}:=\coi_{r}-\coi_{r-1}-1, \quad
\parathr_{r}:=\lambda_{i_{r}}-\lambda_{i_{r+1}}-1, \quad
\parafir_{r}:=\lambda_{i_{r}}+\coi_{r}-\rnk-1
\end{align}
where we write $\coi_0=0$, and let
\begin{align}\label{definition of y}
\parafor_r:=
\binom{ \parasec_{r} }{ \parafir_{r} }
\binom{ \parathr_{r} }{ \parafir_{r} }
\quad \text{for $r=1,\cdots,s$}.
\end{align}
See Figure \ref{Intro picture} for the pictorial meaning of these numbers.
We use the convention $\binom{x}{y}=0$ unless $0\leq y\leq x$.
Now, let 
\begin{align}\label{intro int num}
I(\lambda) := (-1)^{\rnk+s}\parafor_{1}\cdots \parafor_{s}.
\end{align}
The next is the main theorem of this section.

\begin{theorem}\label{intro formula of type A intersection}
If $\A_1,\cdots,\A_{\rnk}$ form a nested chain of subsets, 
then we have
\begin{align*}
\intnum{\tau_{\A_1}\cdots\tau_{\A_{\rnk}}} 
= I(\lambda)
\end{align*}
where $\mu_X$ is the fundamental homology class and $\lambda$ is the Young diagram consisting of $|\A_1|,\cdots,|\A_{\rnk}|$ reordered as a weakly decreasing sequence. Otherwise, the intersection number is zero.
\end{theorem}\vspace{5pt}
\proof
Recall that $\lambda$ is the Young diagram defined by $\lambda_i=|S_{\rnk+1-i}|$ for $i=1,\cdots,\rnk$.
We denote $J(\lambda) := \intnum{\tau_{\A_1}\cdots\tau_{\A_{\rnk}}}$, and we show that $J(\lambda)=I(\lambda)$.
Observe that the condition $0 \leq \parafir_{r} \leq \parathr_{r}$ for all $r=1,\cdots,s$ is equivalent to the condition that each step of the zigzag line of the corners of $\lambda$ crosses the dotted anti-diagonal.
If this condition is not satisfied, then both of $J(\lambda)$ and $I(\lambda)$ are zero. Hence, in the following, we can assume that this condition holds.

We prove the claim by induction on $k:=\rnk-s$.
For the case $k=0$, we have $\lambda_i=|S_{\rnk+1-i}|=\rnk+1-i$ for all $1\leq i\leq \rnk$. So we have $J(\lambda)=1$ by (\ref{general type I=1}).
Since $\parafor_1=\cdots=\parafor_{\rnk}=1$ in this case, we have $I(\lambda)=1$, and the claim follows.
For a general case, there is a lower-right corner (say $r$-th corner from the top) of $\lambda$ whose vertical line has length $\geq2$. 
Then, Lemma \ref{invariance for type A} and Proposition \ref{dotted vanishing type A} combined together show that 
\begin{align}\label{recurrence eq type A}
&J(\lambda)
=
\begin{cases}
-
\begin{pmatrix} \parathr_{r-1} \\ \parafir_{r-1} \end{pmatrix} 
J(\lambda') 
-
\begin{pmatrix} \parathr_{r} \\ \parafir_{r} \end{pmatrix} 
J(\lambda'') 
\qquad
&\text{(if $\lambda\neq\lambda',\lambda''$)}
\vspace{10pt} \\
-
\begin{pmatrix} \parathr_{r-1} \\ \parafir_{r-1} \end{pmatrix} 
J(\lambda') 
&\text{(if $\lambda\neq\lambda', \lambda=\lambda''$)}
\vspace{10pt} \\
-
\begin{pmatrix} \parathr_{r} \\ \parafir_{r} \end{pmatrix} 
J(\lambda'')
&\text{(if $\lambda=\lambda', \lambda\neq\lambda''$)}
\end{cases}
\end{align}
where $\parasec_{r}$, $\parathr_{r}$, and $\parafir_{r}$ are those for $\lambda$, and the Young diagrams $\lambda'$ and $\lambda''$ are given by
\begin{align*}
\lambda'_j = 
\begin{cases}
\rnk+1-j \ &\text{if $j=i_{r-1}+1$}, \\
\lambda_j  &\text{otherwise},
\end{cases}
\qquad 
\lambda''_j = 
\begin{cases}
\rnk+1-j \ &\text{if $j=i_r$}, \\
\lambda_j  &\text{otherwise}.
\end{cases}
\end{align*}
Note that there are no cases that $\lambda=\lambda'=\lambda''$ since our vertical line has length $\geq2$.
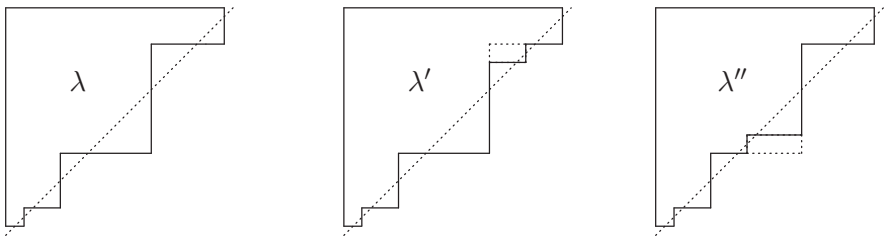
\begin{figure}[h]
\centering
\unitlength 0.1in
\begin{picture}( 45.9500, 11.9300)( 13.7500,-25.6800)
%
\special{pn 8}%
\special{pa 1376 2520}%
\special{pa 1472 2520}%
\special{fp}%
%
\special{pn 8}%
\special{pa 1472 2518}%
\special{pa 1472 2424}%
\special{fp}%
%
\special{pn 8}%
\special{pa 1472 2424}%
\special{pa 1662 2424}%
\special{fp}%
%
\special{pn 8}%
\special{pa 1662 2424}%
\special{pa 1662 2138}%
\special{fp}%
%
\special{pn 8}%
\special{pa 2138 2138}%
\special{pa 2138 1566}%
\special{fp}%
%
\special{pn 8}%
\special{pa 2138 1566}%
\special{pa 2424 1566}%
\special{fp}%
%
\special{pn 8}%
\special{pa 2520 1566}%
\special{pa 2520 1376}%
\special{fp}%
%
\special{pn 8}%
\special{pa 2520 1376}%
\special{pa 1376 1376}%
\special{fp}%
%
\special{pn 8}%
\special{pa 1376 1376}%
\special{pa 1376 2520}%
\special{fp}%
%
\special{pn 8}%
\special{pa 1376 2568}%
\special{pa 2566 1376}%
\special{dt 0.030}%
\put(17.1000,-18.1000){\makebox(0,0)[lb]{$\lambda$}}%
%
\special{pn 8}%
\special{pa 2138 2138}%
\special{pa 1662 2138}%
\special{fp}%
%
\special{pn 8}%
\special{pa 2424 1566}%
\special{pa 2520 1566}%
\special{fp}%
%
\special{pn 8}%
\special{pa 3240 2518}%
\special{pa 3240 2424}%
\special{fp}%
%
\special{pn 8}%
\special{pa 3240 2424}%
\special{pa 3432 2424}%
\special{fp}%
%
\special{pn 8}%
\special{pa 3432 2424}%
\special{pa 3432 2138}%
\special{fp}%
%
\special{pn 8}%
\special{pa 4098 1566}%
\special{pa 4290 1566}%
\special{fp}%
%
\special{pn 8}%
\special{pa 4290 1566}%
\special{pa 4290 1376}%
\special{fp}%
%
\special{pn 8}%
\special{pa 4290 1376}%
\special{pa 3146 1376}%
\special{fp}%
%
\special{pn 8}%
\special{pa 3146 1376}%
\special{pa 3146 2520}%
\special{fp}%
%
\special{pn 8}%
\special{pa 3908 2138}%
\special{pa 3908 1662}%
\special{fp}%
%
\special{pn 8}%
\special{pa 3908 2138}%
\special{pa 3432 2138}%
\special{fp}%
%
\special{pn 8}%
\special{pa 3908 1662}%
\special{pa 4098 1662}%
\special{fp}%
\special{pa 4098 1662}%
\special{pa 4098 1566}%
\special{fp}%
%
\special{pn 8}%
\special{pa 4098 1662}%
\special{pa 3908 1662}%
\special{pa 3908 1566}%
\special{pa 4098 1566}%
\special{pa 4098 1662}%
\special{dt 0.030}%
%
\special{pn 8}%
\special{pa 4874 2518}%
\special{pa 4874 2424}%
\special{fp}%
%
\special{pn 8}%
\special{pa 4874 2424}%
\special{pa 5066 2424}%
\special{fp}%
%
\special{pn 8}%
\special{pa 5066 2424}%
\special{pa 5066 2138}%
\special{fp}%
%
\special{pn 8}%
\special{pa 5542 1566}%
\special{pa 5922 1566}%
\special{fp}%
%
\special{pn 8}%
\special{pa 5922 1566}%
\special{pa 5922 1376}%
\special{fp}%
%
\special{pn 8}%
\special{pa 5542 2042}%
\special{pa 5542 1566}%
\special{fp}%
%
\special{pn 8}%
\special{pa 5256 2138}%
\special{pa 5066 2138}%
\special{fp}%
%
\special{pn 8}%
\special{pa 5256 2138}%
\special{pa 5256 2042}%
\special{fp}%
\special{pa 5256 2042}%
\special{pa 5542 2042}%
\special{fp}%
%
\special{pn 8}%
\special{pa 5542 2042}%
\special{pa 5256 2042}%
\special{pa 5256 2138}%
\special{pa 5542 2138}%
\special{pa 5542 2042}%
\special{dt 0.030}%
\put(34.8000,-18.1000){\makebox(0,0)[lb]{$\lambda'$}}%
\put(51.0000,-18.1000){\makebox(0,0)[lb]{$\lambda''$}}%
%
\special{pn 8}%
\special{pa 3146 2520}%
\special{pa 3240 2520}%
\special{fp}%
%
\special{pn 8}%
\special{pa 5922 1376}%
\special{pa 4778 1376}%
\special{fp}%
%
\special{pn 8}%
\special{pa 4778 1376}%
\special{pa 4778 2520}%
\special{fp}%
%
\special{pn 8}%
\special{pa 4778 2520}%
\special{pa 4874 2520}%
\special{fp}%
%
\special{pn 8}%
\special{pa 3146 2568}%
\special{pa 4336 1376}%
\special{dt 0.030}%
%
\special{pn 8}%
\special{pa 4778 2568}%
\special{pa 5970 1376}%
\special{dt 0.030}%
\end{picture}%
\caption{The Young diagrams $\lambda$, $\lambda'$ and $\lambda''$}
\label{lambda and lambda' and lambda''}
\end{figure}

\noindent
If $\lambda\neq\lambda'$, by the induction hypothesis, we have
\begin{align}\label{inductive comp 10}
J(\lambda') 
&= 
(-1)^{\rnk+(s+1)}
\parafor_1\cdots \parafor_{r-2} \cdot 
\binom{\parasec_{r-1}}{\parafir_{r-1}}
\cdot 1 \cdot
\binom{\parasec_{r}-1}{\parafir_{r}}
\binom{\parathr_{r}}{\parafir_{r}}
\cdot \parafor_{r+1}\cdots \parafor_{s}.
\end{align}
Similarly, if $\lambda\neq\lambda''$, we have
\begin{align}\label{inductive comp 20}
J(\lambda'')
&= 
(-1)^{\rnk+(s+1)}
\parafor_1\cdots \parafor_{r-1} \cdot 
\binom{\parasec_{r}-1}{\parafir_{r}-1}
\cdot 1 \cdot 
\parafor_{r+1}\cdots \parafor_{s}.
\end{align}
Observe that the right-hand-sides of (\ref{inductive comp 10}) and (\ref{inductive comp 20}) vanish when $\lambda=\lambda'$ and $\lambda=\lambda''$, respectively.
Hence, the right-hand-side of the equation (\ref{recurrence eq type A}) can be written as
\begin{align*}
J(\lambda)
&=
-
\binom{\parathr_{r-1}}{\parafir_{r-1}} \cdot
(-1)^{\rnk+s+1}\parafor_1\cdots \parafor_{r-2}
\binom{\parasec_{r-1}}{\parafir_{r-1}}
\binom{\parasec_{r}-1}{\parafir_{r}}
\binom{\parathr_{r}}{\parafir_{r}}
\parafor_{r+1}\cdots \parafor_{s}  \\
&\qquad\qquad
-
\binom{\parathr_{r}}{\parafir_{r}}
\cdot
(-1)^{\rnk+s+1}\parafor_1\cdots \parafor_{r-2}
\parafor_{r-1} \binom{\parasec_{r}-1}{\parafir_{r}-1} 
\parafor_{r+1}\cdots \parafor_{s} \\
&=(-1)^{\rnk+s}\parafor_1\cdots \parafor_{s}=I(\lambda).
\end{align*}
\qed

\vspace{10pt}
\section{The ring structure of the cohomology}\label{ring str of coh}
The cohomology ring $H^*(X)$ of the toric manifold $X$ associated with the fan $\Delta(A_{\rnk})$ is given by the face ring of $\Delta(A_{\rnk})$ modulo the linear relations (\ref{prelim linear relation}) (\cite{Fulton1}). 
As an application of Theorem \ref{intro formula of type A intersection}, we describe the ring structure of the cohomology $H^*(X)$ 
in terms of an additive basis.

Recall that $\invdiv{u\fcw_i}$ for some $i\in[\rnk]$ and permutation $u\in \mathfrak{S}_{\rnk+1}$ is the invariant divisor of $X$ associated with the ray generated by $u\fcw_i\in\CL$. 
Let
\begin{align*}
X_{u} := \bigcap_{i} \invdiv{u\omega_i}
\end{align*}
for each permutation $u\in \mathfrak{S}_{\rnk+1}$ where 
$i$ runs over all descents in $u$. 
Here, a descent in $u$ is a number $i\in[\rnk]$ which satisfies $u(i)>u(i+1)$, and we denote by $d(u)$ the number of descents in $u$.
Denote by $[X_u]\in H^{2d(u)}(X)$ the Poincar$\acute{\text{e}}$ dual of $X_{u}$, then we have
\begin{align*}
[X_{u}]
 = \prod_{i} \tau_{\{u(1),\cdots,u(i)\}}
\end{align*}
where $i$ runs over all descents in $u$ since invariant divisors of $X$ intersect transversely.
$\{[X_u]\}_{u\in \mathfrak{S}_{\rnk+1}}$ forms a module basis of $H^*(X)$ (See \cite{Klyachko} or \cite{Batyrev-Blume} for combinatorial proofs and \cite{De Mari-Procesi-Shayman} for a geometric proof). 
The class $[X_u]$ can be expressed by a Young diagram consisting of the descents in $u$ with the numbers in the nested chain of subsets in 
$D(u)$ (see (\ref{def of D(u)}) for the definition)
written above the diagram so that each column represents the written number above it. 
This expression effectively encodes the descents in $u$ and the information of the chain of subsets.
Denoting $Y^{u}:=w_0X_{w_0u}=\cap X_{u[i]}$ where the intersection runs over all ascents $i$ in $u$, the similar expression works for $[Y^u]$ and the chain of subsets in 
$A(u)$.
Here, $w_0$ is the longest element of $\mathfrak{S}_{\rnk+1}$.
\begin{figure}[h]
\centering
\unitlength 0.1in
\begin{picture}( 33.6000,  6.1000)( 12.0000,-18.7000)
%
\special{pn 8}%
\special{pa 2140 1590}%
\special{pa 2000 1590}%
\special{pa 2000 1450}%
\special{pa 2140 1450}%
\special{pa 2140 1590}%
\special{fp}%
%
\special{pn 8}%
\special{pa 2280 1590}%
\special{pa 2140 1590}%
\special{pa 2140 1450}%
\special{pa 2280 1450}%
\special{pa 2280 1590}%
\special{fp}%
%
\special{pn 8}%
\special{pa 2420 1590}%
\special{pa 2280 1590}%
\special{pa 2280 1450}%
\special{pa 2420 1450}%
\special{pa 2420 1590}%
\special{fp}%
%
\special{pn 8}%
\special{pa 2560 1590}%
\special{pa 2420 1590}%
\special{pa 2420 1450}%
\special{pa 2560 1450}%
\special{pa 2560 1590}%
\special{fp}%
\put(20.3000,-14.3000){\makebox(0,0)[lb]{$2$}}%
\put(21.7000,-14.3000){\makebox(0,0)[lb]{$1$}}%
\put(23.1000,-14.3000){\makebox(0,0)[lb]{$6$}}%
\put(24.5000,-14.3000){\makebox(0,0)[lb]{$4$}}%
%
\special{pn 8}%
\special{pa 2140 1730}%
\special{pa 2000 1730}%
\special{pa 2000 1590}%
\special{pa 2140 1590}%
\special{pa 2140 1730}%
\special{fp}%
%
\special{pn 8}%
\special{pa 4140 1590}%
\special{pa 4000 1590}%
\special{pa 4000 1450}%
\special{pa 4140 1450}%
\special{pa 4140 1590}%
\special{fp}%
%
\special{pn 8}%
\special{pa 4280 1590}%
\special{pa 4140 1590}%
\special{pa 4140 1450}%
\special{pa 4280 1450}%
\special{pa 4280 1590}%
\special{fp}%
%
\special{pn 8}%
\special{pa 4420 1590}%
\special{pa 4280 1590}%
\special{pa 4280 1450}%
\special{pa 4420 1450}%
\special{pa 4420 1590}%
\special{fp}%
\put(12.0000,-15.9000){\makebox(0,0)[lb]{$[X_{216435}]:$}}%
\put(32.0000,-15.9000){\makebox(0,0)[lb]{$[Y^{534162}]:$}}%
%
\special{pn 8}%
\special{pa 4140 1730}%
\special{pa 4000 1730}%
\special{pa 4000 1590}%
\special{pa 4140 1590}%
\special{pa 4140 1730}%
\special{fp}%
%
\special{pn 8}%
\special{pa 2280 1730}%
\special{pa 2140 1730}%
\special{pa 2140 1590}%
\special{pa 2280 1590}%
\special{pa 2280 1730}%
\special{fp}%
%
\special{pn 8}%
\special{pa 2420 1730}%
\special{pa 2280 1730}%
\special{pa 2280 1590}%
\special{pa 2420 1590}%
\special{pa 2420 1730}%
\special{fp}%
%
\special{pn 8}%
\special{pa 2140 1870}%
\special{pa 2000 1870}%
\special{pa 2000 1730}%
\special{pa 2140 1730}%
\special{pa 2140 1870}%
\special{fp}%
\put(40.3000,-14.3000){\makebox(0,0)[lb]{$5$}}%
\put(41.7000,-14.3000){\makebox(0,0)[lb]{$3$}}%
\put(43.1000,-14.3000){\makebox(0,0)[lb]{$4$}}%
%
\special{pn 8}%
\special{pa 4560 1590}%
\special{pa 4420 1590}%
\special{pa 4420 1450}%
\special{pa 4560 1450}%
\special{pa 4560 1590}%
\special{fp}%
\put(44.5000,-14.3000){\makebox(0,0)[lb]{$1$}}%
%
\special{pn 8}%
\special{pa 4280 1730}%
\special{pa 4140 1730}%
\special{pa 4140 1590}%
\special{pa 4280 1590}%
\special{pa 4280 1730}%
\special{fp}%
\end{picture}%
\caption{Two examples for $\rnk=5$ in one-line notations}
\label{X_u and Young diagrams}
\end{figure}
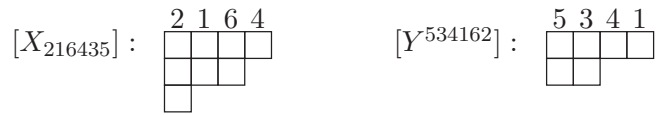

For $u,v,w\in \mathfrak{S}_{\rnk+1}$, we have the Young diagram $\lambda_{uv}^w$ constructed in Section \ref{introduction}.
Recall that $\mu_X$ is the fundamental homology class of $X$.
The following corollary provides the combinatorial rule to compute the intersection number of $X_u$, $X_v$, and $Y^w$ in $X$.
%
\begin{corollary}\label{triple intersection}
$\displaystyle{\intnum{[Y^w][X_u][X_v]}=I(\WY{u}{v}{w})}$.
\end{corollary}
For example, for $\rnk=4$, we have 
\begin{align*}
\intnum{[Y^{35421}][X_{12354}][X_{31254}]} = 2.
\end{align*}
In Figure \ref{lambda and lambda' and lambda''}, we left the numbers on the Young diagram so that we can see the nested chain of subsets appeared in the construction of $\WY{u}{v}{w}$.
\begin{figure}[h]
\centering
\unitlength 0.1in
\begin{picture}( 15.6000,  7.5000)( 16.0000,-18.4000)
%
\special{pn 8}%
\special{sh 0.300}%
\special{pa 3160 1560}%
\special{pa 3020 1560}%
\special{pa 3020 1420}%
\special{pa 3160 1420}%
\special{pa 3160 1560}%
\special{ip}%
%
\special{pn 8}%
\special{sh 0.300}%
\special{pa 2740 1840}%
\special{pa 2600 1840}%
\special{pa 2600 1700}%
\special{pa 2740 1700}%
\special{pa 2740 1840}%
\special{ip}%
%
\special{pn 8}%
\special{pa 2600 1280}%
\special{pa 2740 1280}%
\special{pa 2740 1840}%
\special{pa 2600 1840}%
\special{pa 2600 1280}%
\special{fp}%
%
\special{pn 8}%
\special{pa 2600 1280}%
\special{pa 3160 1280}%
\special{pa 3160 1560}%
\special{pa 2600 1560}%
\special{pa 2600 1280}%
\special{fp}%
%
\special{pn 8}%
\special{pa 3160 1560}%
\special{pa 2600 1560}%
\special{pa 2600 1420}%
\special{pa 3160 1420}%
\special{pa 3160 1560}%
\special{fp}%
%
\special{pn 8}%
\special{pa 2600 1700}%
\special{pa 2740 1700}%
\special{fp}%
%
\special{pn 8}%
\special{pa 2880 1560}%
\special{pa 3020 1560}%
\special{pa 3020 1280}%
\special{pa 2880 1280}%
\special{pa 2880 1560}%
\special{fp}%
%
\special{pn 8}%
\special{pa 2600 1840}%
\special{pa 3160 1280}%
\special{fp}%
\put(26.3000,-12.6000){\makebox(0,0)[lb]{$3$}}%
\put(27.7000,-12.6000){\makebox(0,0)[lb]{$1$}}%
\put(29.1000,-12.6000){\makebox(0,0)[lb]{$2$}}%
\put(30.5000,-12.6000){\makebox(0,0)[lb]{$5$}}%
\put(16.0000,-15.8000){\makebox(0,0)[lb]{$\lambda_{12354, 31254}^{35421}=$}}%
\end{picture}%
\caption{}
\label{lambda and lambda' and lambda''}
\end{figure}

Since $\{[X_u]\}_{u\in \mathfrak{S}_{\rnk+1}}$ forms a module basis of $H^*(X)$, we can consider the expansion coefficients of the product
\begin{align}\label{expansion of product}
[X_{u}][X_{v}] = \sum_{w} \SC{u}{v}{w} [X_{w}].
\end{align}
For example, these coefficients for the product $[X_{s_i}][X_{s_j}]$ can be calculated directly if $|i-j|\geq1$ where $s_i$ is the simple reflection exchanging $i$ and $i+1$.
In fact, we have 
\begin{align}\label{direct computations}
[X_{s_i}][X_{s_j}] = 
\begin{cases}
[X_{s_is_j}](=[X_{s_js_i}]) \quad &\text{if $|i-j|\geq2$}, \\
0 &\text{if $|i-j|=1$}.
\end{cases}
\end{align}
since $\{1,\cdots,i-1,i+1\}$ and $\{1,\cdots,i,i+2\}$ do not form a chain of subsets.

Let us describe each structure constant $c_{u,v}^w$ in terms of intersection numbers computed above. Since a Weyl chamber $\sigma_{u}=\text{cone}(u\fcw_1,\cdots,u\fcw_{\rnk})$ is a maximal cone of the fan, $\sigma_{u}$ is identified with a fixed point of the canonical torus action on $X$ denoted by $p_u\in X$ where $p_u$ is the intersection $\cap_{i=1}^{\rnk}\invdiv{u\omega_i}$.
Then from the definition of $X_{w'}$, one can show that 
$p_u\in X_{w'}$ implies $u\geq w'$ (e.g. \cite{Bjorner-Brenti}; Theorem 2.6.3) where $>$ is the Bruhat order.
If $Y^{w}\cap X_{w'}\neq\emptyset$, then $Y^{w}\cap X_{w'}$ must contain a fixed point since it is an intersection of invariant divisors of $X$, and hence it follows that $w\geq w'$. 
From this observation, we see that $Y^{w}\cap X_{w'}=\emptyset$ unless $w\geq w'$. Also, it is easy to see that $Y^w$ and $X_{w'}$ intersect transversally when $w=w'$. 
Recalling that the class $[Y^{w}][X_{w'}]$ is supported on the intersection $Y^{w}\cap X_{w'}$, we obtain 
\begin{align}\label{intersection of up and down}
\intnum{[Y^{w}][X_{w'}]}=
\begin{cases}
0 \quad &\text{unless $w\geq w'$ and $d(w)=d(w')$} , \\
1 &\text{if } w=w'.
\end{cases} 
\end{align}
See \cite{De Mari-Procesi-Shayman} for a proof using a cellular decomposition of $X$.
Let $\imat$ be the matrix whose $(u,v)$-component is given by $\imatcomp{u}{v}=\intnum{[Y^u][X_v]}=I(\WY{v}{\hspace{1pt}\text{id}}{u})$ for all $u,v\in W$. This matrix $\imat$ is invertible over $\Z$ because of (\ref{intersection of up and down}). 
Now, each coefficient $\SC{u}{v}{w}$ in (\ref{expansion of product}) is a linear transform of the intersection numbers $I(\WY{u}{v}{w})$;
\begin{align}\label{formula for str const}
\SC{u}{v}{w}
= \sum_{w'} \imatcompinv{w}{w'} I(\WY{u}{v}{w'}).
\end{align}
We note that it suffices to take the sum for $w'$ satisfying $d(w)=d(w')$ and $w\geq w'$ since $\imatcompinv{w}{w'}$ is also upper-triangular in the sense of the right-hand-side of (\ref{intersection of up and down}). 

\noindent
%
So the formula (\ref{formula for str const}) exhibits the upper-triangularity of $\SC{u}{v}{w}$ in the sense that $\SC{u}{v}{w}=0$ unless $u,v\leq w$ since $I(\WY{u}{v}{w})$ satisfies the same property.

The transition formula (\ref{formula for str const}) together with (\ref{intersection of up and down}) provides us a recursive formula for the structure constants $\SC{u}{v}{w}$ which is manifestly integral;  
\begin{align*}
\SC{u}{v}{w} = I(\WY{u}{v}{w}) - \sum_{w>w'}\imatcomp{w}{w'}\SC{u}{v}{w'}.
\end{align*} 
Note again that it is enough to take the sum for all $w'$ satisfying $d(w)=d(w')$ and $w> w'$.
From this recursion, we recover (\ref{direct computations}), and we can compute the expansion of $[X_{s_i}]^2$. For example, if $\rnk=3$, we obtain
\begin{align*}
[X_{2134}]^2
&= [X_{2431}] - [X_{4213}] - [X_{3421}] - [X_{3241}] - [X_{3214}]. 
\end{align*}

\section{Other classical types}
Note that the argument in the previous section can be naturally generalized to arbitrary root systems by considering the non-singular subvariety
\begin{align}\label{def of Xu}
X_u = \bigcap_{i} \invdiv{u\fcw_i}
\end{align} 
for each $u\in W$ where $\invdiv{u\fcw_i}$ is the invariant divisor of $X$ corresponding to the ray generated by $u\fcw_i$ and $i$ runs over all $i$ satisfying $u(\alpha_i)\in\Phi^-$.
Here, $\Phi^-$ is the set of negative roots. It follows that the Poincar$\acute{\text{e}}$ duals $\{[X_u]\}_{u\in W}$ form an additive basis of the integral cohomlogy $H^*(X)$ (see \cite{De Mari-Procesi-Shayman}).
\begin{remark}
\emph{
The collections $\{c_{u, v}^w\}_{u,v,w\in W}$ and $\{\intnum{[Y^w][X_u][X_v]}\}_{u,v,w\in W}$ are independent on the choice of the simple roots $\Pi$.
}
\end{remark}

\subsection{Intersection numbers for type $B_{\rnk}$}\label{section int num type B}
For the classical root system of type $B_{\rnk}$, the roots are $\{t_i-t_j, \ \pm(t_i+t_j) , \ \pm t_{i} \in E \mid 1\leq i\neq j\leq \rnk\}$ where $E=\R^{\rnk}$. We choose $\SimR=\{t_i-t_{i+1}, \ t_{\rnk} \mid 1\leq i\leq \rnk-1\}$ as a set of simple roots, and write $\alpha_i=t_i-t_{i+1} (1\leq i\leq\rnk-1)$, $\alpha_{\rnk}=t_{\rnk}$. The Weyl group $\widetilde{\mathfrak{S}}_{\rnk}$ is the $\rnk$-th signed permutation group. Letting $t_{-i}:=-t_i$ for all $1\leq i\leq \rnk$, $u\in\widetilde{\mathfrak{S}}_{\rnk}$ acts on $E$ by $ut_i=t_{u(i)}$.
The minimal generators $\fcw_1,\cdots,\fcw_{\rnk}\in E^*$ of the fundamental Weyl chamber are $\fcw_i=e_1+\cdots+e_i$ for $i=1,\cdots,\rnk$.

Let $[\pm\rnk] = \{1,\cdots,\rnk, -1,\cdots,-\rnk\}$.
For $\A\in2^{[\pm\rnk]}$, consider a condition
\begin{align}\label{type B basics 100}\tag{$*$}
\text{ for any $i\in[\pm\rnk]$, if $i\in \A$ then $-i\notin \A$}.
\end{align} 
We have a well-defined map $\dPhi \rightarrow 2^{[\pm\rnk]}$ by $u\fcw_i \mapsto \{u(1),\cdots,u(i)\}$.
This leads us to an identification
\begin{align*}
\dPhi \quad \longleftrightarrow \quad \text{the set of non-empty subsets of $[\pm\rnk]$ satisfying (\ref{type B basics 100})}.
\end{align*} 
Now, for each $\emptyset\subsetneq S\subset[\pm\rnk]$ satisfying (\ref{type B basics 100}), we define
$\tau_{\A}:=\tau_{u\omega_i}$ where $u\omega_i\in\dPhi$ corresponds to $\A$ by this identification. 
For $\emptyset\subsetneq \A_1,\cdots, \A_{q}\subset[\pm\rnk] \ (1\leq q\leq\rnk)$ satisfying (\ref{type B basics 100}), we have that $\tau_{\A_1}\cdots\tau_{\A_{q}}=0$ unless these sets form a nested chain of subsets, as in the case for type $A_{\rnk}$.

For each $k\in[\pm\rnk]$, let $\Bsetsec\subset[\pm\rnk]$ satisfy (\ref{type B basics 100}), $k\in \Bsetsec$, and $|\Bsetsec|=\rnk$. From the linear relation (\ref{prelim linear relation}) for the root $\alpha=t_k$, we can deduce that
\begin{align*}
{\tau_{\Bsetsec}}^2
= -
\sum_{\substack{k\in \Bsetrun \\ \Bsetrun\subsetneq \Bsetsec}} \tau_{\Bsetrun} \tau_{\Bsetsec}
\end{align*}  
where the sum is taken over all $\emptyset \subsetneq \Bsetrun\subset [\pm\rnk]$ satisfying (\ref{type B basics 100}) with the prescribed conditions.
Similarly, for each $k,l\in[\pm\rnk]$, let $\Bsetsec\subset[\pm\rnk]$ satisfy (\ref{type B basics 100}), $k\in \Bsetsec$, and $\pm l\notin \Bsetsec$ (hence $1\leq |\Bsetsec|\leq\rnk-1$). 
Then from (\ref{prelim linear relation}) for the root $\alpha=t_k-t_l$, we obtain 
\begin{align*}
&{\tau_{\Bsetsec}}^2
=
-
\sum_{\substack{k\in \Bsetrun, \ \pm l\notin \Bsetrun \\ \Bsetrun\neq \Bsetsec}} \tau_{\Bsetrun} \tau_{\Bsetsec} 
-
\sum_{k,-l\in \Bsetrun} 2\tau_{\Bsetrun}\tau_{\Bsetsec} .
\end{align*} 
Observe that the second summand will vanish after multiplying $\tau_{\Bsetfir}$ and $\tau_{\Bsetthr}$ for $\Bsetfir\subset \Bsetsec\backslash{\{k\}}$ and $\Bsetsec\coprod\{l\}\subset \Bsetthr$ where we write $\tau_{\emptyset}=0$.
So these two equations can be used to prove the separation rule similar to Lemma \ref{type A separation}, and we obtain the same type of vanishing property as in Proposition \ref{dotted vanishing type A}. 
Now the argument in the proof of Theorem \ref{intro formula of type A intersection} also works for this case, and it follows that
\begin{theorem}\label{formula of type B intersection}
If $\emptyset\subsetneq \A_1,\cdots,\A_{\rnk}\subset[\pm\rnk]$ satisfying \emph{(\ref{type B basics 100})} form a nested chain of subsets, then we have
\begin{align*}
\intnum{\tau_{\A_1}\cdots\tau_{\A_{\rnk}}}
= 2^{\rnk-\lambda_1} I(\lambda)
\end{align*}
where 
$\mu_X$ is the fundamental homology class of $X$ and
$\lambda$ is the Young diagram consisting of $|\A_1|,\cdots,|\A_{\rnk}|$ reordered as a weakly decreasing sequence and $I$ is the function defined in \eqref{intro int num}. Otherwise, the intersection number is zero.
\end{theorem}

Let $\alpha_i:=t_{i}-t_{i+1}$ for $1\leq i\leq\rnk-1$ and $\alpha_{\rnk}:=t_{\rnk}$.
For each signed permutation $u\in \widetilde{\mathfrak{S}}_{\rnk}$, an element $i\in[\rnk]$ satisfies $u(\alpha_i)\in\Phi^-$ if and only if
\begin{itemize}
 \item[(D-1)] if $i\leq n-1$, then $u(i)>u(i+1)$ with the same sign or $u(i)<u(i+1)$ with different signs,
 \item[(D-2)] if $i=\rnk$, then $u(i)<0$.
\end{itemize}
Similarly, consider the conditions
\begin{itemize}
 \item[(A-1)] if $i\leq n-1$, then $u(i)<u(i+1)$ with the same sign or $u(i)>u(i+1)$ with different signs,
 \item[(A-2)] if $i=\rnk$, then $u(i)>0$.
\end{itemize}
Denoting 
\begin{align*}
D(u):=\{ u[i] \mid \text{$i$ satisfies (D)} \}  
\quad \text{and} \quad
A(u):=\{ u[i] \mid \text{$i$ satisfies (A)} \},
\end{align*}
we define a Young diagram $\lambda_{u,v}^w$ in the manner described in the last section.
Note that we put $I(\emptyset)=0$ as a convention.

Now, for signed permutations $u,v,w\in \widetilde{\mathfrak{S}}_{\rnk}$, the intersection number of $Y^w$, $X_u$, and $X_v$ in $X$ of type $B_{\rnk}$ is given by the following. 
\begin{corollary}\label{triple intersection for type B}
For signed permutations $u,v,w\in \widetilde{\mathfrak{S}}_{\rnk}$, we have
\[ \intnum{[Y^w][X_u][X_v]}=2^{\rnk-(\lambda_{u,v}^w)_1}I(\WY{u}{v}{w}) \]
where $I$ is the function defined in \eqref{intro int num}. 
\end{corollary}

For example, for $\rnk=4$ with the convention $\bar{k}=-k$, Corollary \ref{triple intersection for type B} computes
\begin{align*}
\intnum{[Y^{2\bar{3}\bar{1}\bar{4}}] [X_{2\bar{3}14}] [X_{2\bar{3}14}]}  = 4.
\end{align*}
\vspace{-20pt}
\begin{figure}[h]
\centering
\unitlength 0.1in
\begin{picture}( 15.9000,  7.5000)( 14.3000,-18.4000)
\put(26.3000,-12.6000){\makebox(0,0)[lb]{$\bar{3}$}}%
%
\special{pn 8}%
\special{pa 2740 1560}%
\special{pa 2600 1560}%
\special{pa 2600 1420}%
\special{pa 2740 1420}%
\special{pa 2740 1560}%
\special{fp}%
%
\special{pn 8}%
\special{pa 2740 1420}%
\special{pa 2600 1420}%
\special{pa 2600 1280}%
\special{pa 2740 1280}%
\special{pa 2740 1420}%
\special{fp}%
%
\special{pn 8}%
\special{pa 2600 1560}%
\special{pa 2460 1560}%
\special{pa 2460 1420}%
\special{pa 2600 1420}%
\special{pa 2600 1560}%
\special{fp}%
%
\special{pn 8}%
\special{pa 2600 1420}%
\special{pa 2460 1420}%
\special{pa 2460 1280}%
\special{pa 2600 1280}%
\special{pa 2600 1420}%
\special{fp}%
%
\special{pn 8}%
\special{pa 2600 1700}%
\special{pa 2460 1700}%
\special{pa 2460 1560}%
\special{pa 2600 1560}%
\special{pa 2600 1700}%
\special{fp}%
%
\special{pn 8}%
\special{sh 0.300}%
\special{pa 2600 1840}%
\special{pa 2460 1840}%
\special{pa 2460 1700}%
\special{pa 2600 1700}%
\special{pa 2600 1840}%
\special{fp}%
%
\special{pn 8}%
\special{sh 0.300}%
\special{pa 2740 1700}%
\special{pa 2600 1700}%
\special{pa 2600 1560}%
\special{pa 2740 1560}%
\special{pa 2740 1700}%
\special{fp}%
%
\special{pn 8}%
\special{pa 2460 1840}%
\special{pa 3020 1280}%
\special{fp}%
\put(24.9000,-12.6000){\makebox(0,0)[lb]{$2$}}%
\put(14.3000,-16.6000){\makebox(0,0)[lb]{$\lambda_{2\bar{3}14, 2\bar{3}14}^{2\bar{3}\bar{1}\bar{4}}\ = $}}%
\end{picture}%
\caption{}
\label{EX for type B}
\end{figure}

\subsection{Intersection numbers for type $C_{\rnk}$}\label{section int num type C}
For the classical root system of type $C_{\rnk}$, the roots are $\{t_i-t_j, \ \pm(t_i+t_j), \ \pm2t_i \in E \mid 1\leq i\neq j\leq \rnk\}$ where $E=\R^{\rnk}$. We choose $\SimR=\{t_i-t_{i+1}, \ 2t_{\rnk} \mid 1\leq i\leq \rnk-1\}$ as a set of simple roots, and write $\alpha_i=t_i-t_{i+1} (1\leq i\leq\rnk-1)$, $\alpha_{\rnk}=2t_{\rnk}$. The Weyl group $\widetilde{\mathfrak{S}}_{\rnk}$ is the $\rnk$-th signed permutation group as above.
The minimal generators $\fcw_1,\cdots,\fcw_{\rnk}$ of the fundamental Weyl chamber are $\fcw_i=e_1+\cdots+e_i$ for $i=1,\cdots,\rnk-1$ and $\fcw_{\rnk}=\frac{1}{2}(e_1+\cdots+e_{\rnk})$.

We have a well-defined map $\dPhi \rightarrow 2^{[\pm\rnk]}$ by $v\fcw_i \mapsto \{v(1),\cdots,v(i)\}$, and obtain an identification $\dPhi$ and the set of non-empty subsets of $[\pm\rnk]$ satisfying (\ref{type B basics 100}).
For $\emptyset\subsetneq \A_1,\cdots, \A_{q}\subsetneq[\pm\rnk] \ (1\leq q\leq \rnk)$ satisfying (\ref{type B basics 100}), we have that $\tau_{\A_1}\cdots\tau_{\A_{q}}=0$ unless these sets form a nested chain of subsets where $\tau_{\A}$ is defined as in Section \ref{section int num type B}.

For each $k\in[\pm\rnk]$, let $\Csetsec\subset[\pm\rnk]$ satisfy (\ref{type B basics 100}), $k\in \Csetsec$, and $|\Csetsec|=\rnk$.
Then (\ref{prelim linear relation}) for the root $\alpha=2t_k$ shows that
\begin{align*}
{\tau_{\Csetsec}}^2
= -
\sum_{\substack{k\in \Csetrun \\ \Csetrun\subsetneq \Csetsec}} 2\tau_{\Csetrun} \tau_{\Csetsec}
\end{align*}  
where the sum is taken over all $\emptyset \subsetneq \Bsetrun\subset [\pm\rnk]$ satisfying (\ref{type B basics 100}) with the prescribed conditions.
For each $k,l\in[\pm\rnk]$, let $\Csetsec\subset[\pm\rnk]$ satisfy (\ref{type B basics 100}), $k\in \Csetsec$, and $\pm l\notin \Csetsec$ (hence $1\leq |\Csetsec|\leq\rnk-1$). 
Then from (\ref{prelim linear relation}) for the root $\alpha=t_k-t_l$, we obtain
\begin{align*}
&{\tau_{\Csetsec}}^2
=
-
\sum_{\substack{k\in \Csetrun, \ \pm l\notin \Csetrun \\ \Csetrun\neq \Csetsec}} \tau_{\Csetrun} \tau_{\Csetsec} 
-
\sum_{\substack{k\in \Csetrun, \ -l\in \Csetrun \\ |\Csetrun|\neq\rnk}} 2\tau_{\Csetrun}\tau_{\Csetsec} 
-
\sum_{\substack{k\in \Csetrun, \ -l\in \Csetrun \\ |\Csetrun|=\rnk}} \tau_{\Csetrun} \tau_{\Csetsec}.
\end{align*} 
With a similar observation made for type $B_{\rnk}$, we again have the same type of vanishing property as in Proposition \ref{dotted vanishing type A}. Hence, we obtain 
\begin{theorem}\label{formula of type B intersection}
If $\emptyset\subsetneq \A_1,\cdots,\A_{\rnk}\subsetneq[\pm\rnk]$ satisfying \emph{(\ref{type B basics 100})} form a nested chain of subsets, then we have
\begin{align*}
\intnum{\tau_{\A_1}\cdots\tau_{\A_{\rnk}}}
= 2^{\rnk-\lambda_1+m-1} I(\lambda)
\end{align*}
where $\mu_X$ is the fundamental homology class of $X$ and $\lambda$ is the Young diagram consisting of the numbers $|\A_1|,\cdots,|\A_{\rnk}|$ reordered as a weakly decreasing sequence and $I$ is the function defined in \eqref{intro int num} and $m$ is the number of rows of $\lambda$ of length $\rnk$. Otherwise, the intersection number is zero.
\end{theorem}

For signed permutations $u,v,w\in \widetilde{\mathfrak{S}}_{\rnk}$, let $\lambda_{u,v}^w$ be the Young diagram defined in Section \ref{section int num type B}. 
The intersection number of $Y^w$, $X_u$, and $X_v$ in $X$ of type $C_{\rnk}$ is given by the following.  
\begin{corollary}\label{triple intersection for type C}
For signed permutations $u,v,w\in \widetilde{\mathfrak{S}}_{\rnk}$, we have
\[ \intnum{[Y^w][X_u][X_v]}=2^{\rnk-(\lambda_{u,v}^w)_1+m-1}I(\WY{u}{v}{w}) \]
where $I$ is the function defined in \eqref{intro int num} and and $m$ is the number of rows of $\lambda_{u,v}^w$ of length $\rnk$.
\end{corollary}

\vspace{10pt}
\subsection{Intersection numbers for type $D_{\rnk}$}
For the classical root system of type $D_{\rnk}$, the roots are $\{t_i-t_j, \ \pm(t_i+t_j) \in E \mid 1\leq i\neq j\leq \rnk\}$ where $E=\R^{\rnk}$. We choose $\SimR=\{t_i-t_{i+1}, \ t_{\rnk-1}+t_{\rnk} \mid 1\leq i\leq \rnk-1\}$ as a set of simple roots, and write $\alpha_i=t_i-t_{i+1} (1\leq i\leq\rnk-2)$, $\alpha_{\rnk-1}=t_{\rnk-1}+t_{\rnk}$, $\alpha_{\rnk}=t_{\rnk-1}-t_{\rnk}$. The Weyl group $\widetilde{\mathfrak{S}}_{\rnk}^+$ is the $\rnk$-th \textit{even signed permutation group} defined by
\begin{align*}
\widetilde{\mathfrak{S}}_{\rnk}^+:=\{w\in\widetilde{\mathfrak{S}}_{\rnk} \mid \text{ the number of $i$ with $w(i)<0$ is even}\}
\end{align*}
where $\widetilde{\mathfrak{S}}_{\rnk}$ is the $\rnk$-th signed permutation group.
The minimal generators $\fcw_1,\cdots,\fcw_{\rnk}$ $\in E^*$ of the fundamental Weyl chamber are $\fcw_i=e_1+\cdots+e_i$ for $i=1,\cdots,\rnk-2$, \  $\fcw_{\rnk-1}=\frac{1}{2}(e_1+\cdots+e_{\rnk-1}+e_{\rnk})$ and $\fcw_{\rnk}=\frac{1}{2}(e_1+\cdots+e_{\rnk-1}-e_{\rnk})$.

For $\A\in2^{[\pm\rnk]}$, consider a condition
\begin{align}\label{type D basics 100}\tag{$**$}
\text{$|\A|\neq \rnk-1$, and if $i\in \A$ then $-i\notin \A$ for any $i\in[\pm\rnk]$ }.
\end{align} 
We have a well-defined map $\dPhi \rightarrow 2^{[\pm\rnk]}$ given by
\begin{align*}
&u\fcw_i \mapsto u[i]=\{u(1),\cdots,u(i)\} \quad \text{for } 1\leq i\leq \rnk-2, \\
&u\fcw_{\rnk-1} \mapsto u[n]_+=\{u(1),\cdots,u(\rnk-1), u(\rnk)\}, \\
&u\fcw_{\rnk} \mapsto u[n]_-=\{u(1),\cdots,u(\rnk-1),-u(\rnk)\}.
\end{align*} 
where $[\rnk]_+ = \{1,2,\cdots,n-1,n\}$ and $[\rnk]_- = \{1,2,\cdots,n-1,-n\}$.
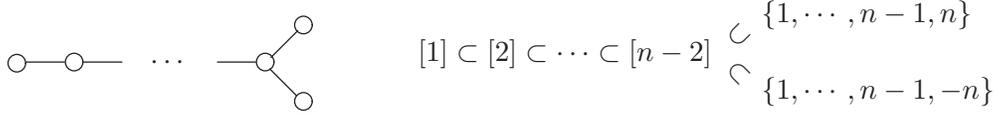
\begin{figure}[h]
\centering
\unitlength 0.1in
\begin{picture}( 54.0000,  6.2000)(  5.3000,-10.8000)
%
\special{pn 8}%
\special{ar 656 806 46 46  0.0000000 6.2831853}%
%
\special{pn 8}%
\special{pa 706 800}%
\special{pa 906 800}%
\special{fp}%
%
\special{pn 8}%
\special{ar 956 800 46 46  0.0000000 6.2831853}%
%
\special{pn 8}%
\special{pa 1006 800}%
\special{pa 1206 800}%
\special{fp}%
%
\special{pn 8}%
\special{pa 1706 800}%
\special{pa 1906 800}%
\special{fp}%
%
\special{pn 8}%
\special{ar 1956 800 46 46  0.0000000 6.2831853}%
%
\special{pn 8}%
\special{ar 2156 606 46 46  0.0000000 6.2831853}%
%
\special{pn 8}%
\special{ar 2156 1006 46 46  0.0000000 6.2831853}%
%
\special{pn 8}%
\special{pa 2126 636}%
\special{pa 1986 776}%
\special{fp}%
%
\special{pn 8}%
\special{pa 2120 976}%
\special{pa 1980 836}%
\special{fp}%
\put(13.5500,-8.3500){\makebox(0,0)[lb]{$\cdots$}}%
\put(27.5500,-8.3000){\makebox(0,0)[lb]{$[1] \subset [2] \subset \cdots \subset [\rnk-2]$}}%
\put(45.5500,-6.3000){\makebox(0,0)[lb]{$\{1,\cdots,n-1,\rnk\}$}}%
\put(45.5500,-10.3000){\makebox(0,0)[lb]{$\{1,\cdots,n-1,-\rnk\}$}}%
%
\special{pn 8}%
\special{pa 530 550}%
\special{pa 5930 550}%
\special{pa 5930 1080}%
\special{pa 530 1080}%
\special{pa 530 550}%
\special{ip}%
\put(43.6000,-9.3000){\makebox(0,0)[lb]{$\rotatebox{-45}{$\subset$}$}}%
\put(43.6000,-7.3000){\makebox(0,0)[lb]{$\rotatebox{45}{$\subset$}$}}%
\end{picture}%
\caption{The Dynkin diagram and a basic chain of subsets for type $D_{\rnk}$}
\label{type D chain}
\end{figure}
\vspace{-3pt}
\\ 
It follows that this map $\dPhi \rightarrow 2^{[\pm\rnk]}$ is an injection.
In fact, we cannot have $\{u(1),\cdots,u(\rnk-1), u(\rnk)\}=\{v(1),\cdots,v(\rnk-1),-v(\rnk)\}$ for any $u,v\in\widetilde{\mathfrak{S}}_{\rnk}^+$ since the number of negative integers in the left hand side and the right-hand-side are different, and so
$u\fcw_{\rnk-1}$ and $v\fcw_{\rnk}$ are never mapped to the same element.
The other cases are left to the reader.
So we can make an identification
\begin{align}\label{type D identification}
\dPhi \quad \longleftrightarrow \quad \text{the set of non-empty subsets of $[\pm\rnk]$ satisfying (\ref{type D basics 100})}.
\end{align} 
Hence,  for each $\emptyset\subsetneq \A\subset[\pm\rnk]$ satisfying (\ref{type D basics 100}), we define
$\tau_{\A}:=\tau_{u\omega_i}$ where $u\omega_i\in\dPhi$ corresponds to $\A$ by this identification.

Let us denote by $\mathcal{C}$ the set of chains of subsets $\{\DA_i\}_i$ of $[\pm n]$ of the following form:
there exists $u\in\widetilde{\mathfrak{S}}_{\rnk}^+$ such that $\DA_i=u[i]$ for $1\leq i\leq n-2$, $\DA_{n-1}=u[n]_+$, and $\DA_{n}=u[n]_-$.
Note that $\{\DA_i\}_i$ does not have a set of order $n-1$ and satisfies the same inclusion relation shown in Figure \ref{type D chain}.
A \textit{subchain} $\{\A_i\}_{i}$ of a chain $\{\DA_i\}_{i}$ in $\mathcal{C}$ is a sequence satisfying $\A_j\in\{\DA_i\}_{i}$ for $1\leq j\leq n$ and $S_{j}\subset S_{j'}$ for $1\leq j\leq j'\leq \rnk$ unless $|S_{j}|=|S_{j'}|=\rnk$.
For $\emptyset\subsetneq\A_1,\cdots,\A_{q}\subsetneq[\pm\rnk] \ (1\leq q\leq \rnk)$ satisfying (\ref{type D basics 100}), we have $\tau_{\A_1}\cdots\tau_{\A_{q}}=0$ unless the sequence forms a subchain of a chain in $\mathcal{C}$ up to reordering. 
Let $\{\A_i\}_i$ be a subchain of a chain in $\mathcal{C}$. 
For $\A_i$ satisfying $|\A_i|=\rnk$, we say that $\A_i$ is \textit{even} (resp. \textit{odd}) if the number of negative elements of $\A_i$ is even (resp. odd).
Recall that $\mu_X$ is the fundamental homology class of $X$.
The following is Lemma \ref{prelim Weyl inv} for type $D_{\rnk}$.
\begin{lemma}\label{invariance property for type D}
Let $\{\A_i\}_i$ and $\{\A'_i\}_i$ be subchains of some chains in $\mathcal{C}$. 
If $|\A_i|=|\A'_i|$ for $i=1,\cdots,\rnk$ and the number of even $\A_i$'s and the number of even $\A'_i$'s are the same, then $\intnum{\tau_{\A_1}\cdots\tau_{\A_{\rnk}}} = \intnum{\tau_{\A'_1}\cdots\tau_{\A'_{\rnk}}}$.
\end{lemma}

Let $\{\A_i\}_i$ be a subchain of a chain in $\mathcal{C}$. We denote by $\lambda$ the \textit{signed Young diagram} consisting of $\lambda_i=|\A_{\rnk+1-i}|$ for $i=1,\cdots,\rnk$ where the label of $\lambda$ is defined as follows:
if we have $\lambda_i=\rnk$, then we label this row by $+$ (resp. $-$) if $\A_{\rnk+1-i}$ is even (resp. odd).
Recall from Section \ref{section type A} that the \textit{dotted anti-diagonal line} drawn on the Young diagram is the dotted line shifted down half the length of a single box from the standard anti-diagonal. 
Our first aim is to prove the following.
\begin{proposition}\label{zigzag lemma type D}
\emph{(The vanishing property)}
$\intnum{\tau_{\A_1}\cdots\tau_{\A_{\rnk}}}=0$ unless each step of the zigzag line of the lower-right corners of $\lambda$ crosses the dotted anti-diagonal.
\end{proposition}
\begin{figure}[h]
\centering
\input{pic-typeD-1}
\caption{}
\label{zigzag for type D}
\end{figure}

For each $k,l\in[\pm\rnk]$, let $\Dsetsec\subset[\pm\rnk]$ satisfy (\ref{type D basics 100}), $k\in \Dsetsec$, and $\pm l\notin \Dsetsec$ (hence $|\Dsetsec|\leq\rnk-2$).
By the linear relation (\ref{prelim linear relation}) for the root $\alpha=t_k-t_l$, it follows that
\begin{align}\label{type D vanishing 130}
{\tau_{\Dsetsec}}^2 = 
-
\sum_{\substack{k\in \Dsetrun, \ \pm l\notin \Dsetrun \\ \Dsetrun\neq \Dsetsec}} \tau_{\Dsetrun} \tau_{\Dsetsec}
-
\sum_{\substack{k\in \Dsetrun, \ -l\in \Dsetrun \\ |\Dsetrun|\neq\rnk}} 2\tau_{\Dsetrun} \tau_{\Dsetsec}
-
\sum_{\substack{k\in \Dsetrun, \ -l\in \Dsetrun \\ |\Dsetrun|=\rnk}} \tau_{\Dsetrun} \tau_{\Dsetsec}
\end{align}
where the sum is taken over all $\emptyset \subsetneq \Bsetrun\subset [\pm\rnk]$ satisfying (\ref{type D basics 100}) with the prescribed conditions.
If $\Dsetfir\subsetneq \Dsetsec$ with $k\notin \Dsetfir$, and if $\Dsetsec\subsetneq \Dsetthr$ with 
$l\in \Dsetthr$ then, 
\begin{align}\label{type D vanishing 140}
\tau_{\Dsetfir} {\tau_{\Dsetsec}}^2 \tau_{\Dsetthr} = 
-
\sum_{\substack{k\in \Dsetrun, \ \pm l\notin \Dsetrun \\ \Dsetfir\subsetneq \Dsetrun \subsetneq \Dsetthr, \ \Dsetrun\neq \Dsetsec}} \tau_{\Dsetfir} \tau_{\Dsetrun} \tau_{\Dsetsec} \tau_{\Dsetthr}
\ \ - \ \ 
\delta_{|\Dsetthr|,\rnk}\tau_{\Dsetfir} \tau_\Dsetsec \tau_{(l,-l)\Dsetthr} \tau_{\Dsetthr}.
\end{align} 
where $\delta_{|\Dsetthr|,\rnk}$ is the Kronecker delta.
If $|\Dsetthr|=\rnk$, then after multiplying (\ref{type D vanishing 140}) by $\tau_{\overline{\Dsetthr}}$ where $\overline{\Dsetthr}=(p,-p)\Dsetthr$ for some $p\in \Dsetthr\backslash \{l\}$, we obtain 
\begin{align}\label{type D vanishing 141}
\tau_{\Dsetfir} {\tau_{\Dsetsec}}^2 \tau_{\Dsetthr}\tau_{\overline{\Dsetthr}} = 
-
\sum_{\substack{k\in \Dsetrun, \ \pm l, \pm p\notin \Dsetrun \\ \Dsetfir\subsetneq \Dsetrun \subsetneq \Dsetthr, \ \Dsetrun\neq \Dsetsec}} \tau_{\Dsetfir} \tau_{\Dsetrun} \tau_{\Dsetsec} \tau_{\Dsetthr}\tau_{\overline{\Dsetthr}}.
\end{align} 

Let $\lambda$ as above. 
We denote 
\begin{align*}
&m_+(\lambda):= |\{i \mid \text{$\lambda_i=\rnk$ and the label of $\lambda_i$ is $+$} \}|, \\
&m_-(\lambda):= |\{i \mid \text{$\lambda_i=\rnk$ and the label of $\lambda_i$ is $-$} \}|.
\end{align*}

\begin{lemma}\label{type D vanishing 300}
Suppose that one of the following holds:
\begin{itemize}
 \item[(i)] $m_+(\lambda)=m_-(\lambda)=1$
 \item[(ii)] $(m_+(\lambda),m_-(\lambda))$ is equal to $(1,0)$ or $(0,1)$,
  \item[(iii)] $m_+(\lambda)=m_-(\lambda)=0$.
\end{itemize}
Then $\intnum{\tau_{S_1}\cdots\tau_{S_{\rnk}}}=0$ unless each step of the zigzag line of the corners of $\lambda$ crosses the dotted anti-diagonal.
\end{lemma}
\proof
%
The claim for the case (i) can be proved by (\ref{type D vanishing 140}) and (\ref{type D vanishing 141}) as in the proof of Proposition \ref{dotted vanishing type A}.
For the case (ii), the same argument works together with (\ref{type D vanishing 140}), since we have already proved the claim for the case (i). 
Now, the case (iii) is shown again by the same proof used for Proposition \ref{dotted vanishing type A} together with (\ref{type D vanishing 130}), (\ref{type D vanishing 140}), and the claim for the case (ii).
\qed

\vspace{20pt}
For each $k,l\in[\pm\rnk]$, let $\emptyset\subsetneq \Dsetsec\subsetneq[\pm\rnk]$ satisfy (\ref{type D basics 100}),  $k,l\in \Dsetsec$, and $|\Dsetsec|=\rnk$. 
If $\Dsetfir\subset \Dsetsec\backslash\{k,l\}$, then from the linear relation (\ref{prelim linear relation}) for the root $\alpha=t_k+t_l$, we obtain
\begin{align}\notag
\tau_{\Dsetfir} {\tau_{\Dsetsec}}^2 = 
&
-
\sum_{\substack{k\in \Dsetrun, \ \pm l\notin \Dsetrun}} \tau_{\Dsetfir} \tau_{\Dsetrun} \tau_{\Dsetsec}
-
\sum_{\substack{\pm k\notin \Dsetrun, \ l\in \Dsetrun}} \tau_{\Dsetfir} \tau_{\Dsetrun} \tau_{\Dsetsec}
\\ \label{type D vanishing 900}
&\quad\quad\quad\quad
-
\sum_{\substack{k, l\in \Dsetrun, \\ |\Dsetrun|\neq\rnk}} 2\tau_{\Dsetfir} \tau_{\Dsetrun} \tau_{\Dsetsec}
-
\sum_{\substack{k,l\in \Dsetrun, \ \Dsetrun\neq \Dsetsec, \\ |\Dsetrun|=\rnk}} \tau_{\Dsetfir} \tau_{\Dsetrun} \tau_{\Dsetsec}
\end{align} 
where we denote $\tau_{\emptyset}=1$.
Especially if $\Dsetfir= \Dsetsec\backslash\{k,l\}$, then  $|\Dsetfir|=\rnk-2$ and we have
\begin{align}\label{type D vanishing 1000}
\tau_{\Dsetfir}{\tau_{\Dsetsec}}^2
&= 
-
\sum_{\substack{k,l\in \Dsetrun, \ \Dsetrun\neq \Dsetsec, \\ |\Dsetrun|=\rnk}} \tau_{\Dsetfir} \tau_{\Dsetrun} \tau_{\Dsetsec} =0.
\end{align} 
The second equality follows since an element of $\Dsetrun$ which is neither $k$ nor $l$ has to be $-1$ times an element of $\Dsetsec$, which implies that $\Dsetfir\not\subset \Dsetrun$. 
On the other hand, letting 
$\overline{\Dsetsec}=(-k,k)\Dsetsec$, we obtain from (\ref{type D vanishing 900}) that 
\begin{align}\label{type D vanishing 1100}
\tau_{\Dsetfir} {\tau_\Dsetsec}^2 \tau_{\overline{\Dsetsec}} = 
&
-
\sum_{\substack{\pm k\notin \Dsetrun, \ l\in \Dsetrun \\ \Dsetfir\subsetneq \Dsetrun\subsetneq \Dsetsec}} \tau_{\Dsetfir} \tau_{\Dsetrun} \tau_{\Dsetsec} \tau_{\overline{\Dsetsec}}.
\end{align}

\vspace{5pt}
\begin{lemma}\label{type D vanishing 1200}
Suppose that one of the following holds:
\begin{itemize}
 \item[(i)] $m_+(\lambda), m_-(\lambda)\geq1$,
 \item[(ii)] $m_+(\lambda)\geq1$ and $m_-(\lambda)=0$, 
 \item[(iii)] $m_+(\lambda)=0$ and $m_-(\lambda)\geq1$.
\end{itemize}
Then $\intnum{\tau_{S_1}\cdots\tau_{S_{\rnk}}}=0$ unless each step of the zigzag line of the lower-right corners of $\lambda$ crosses the dotted anti-diagonal.
\end{lemma}
\proof
The claim for the case (i) follows from (\ref{type D vanishing 1100}) and the case (i) of Lemma \ref{type D vanishing 300} by induction on $m_+(\lambda)+m_-(\lambda)$.
Let us consider the case (ii). 
We prove the claim by induction on $m_+(\lambda)$.
For the case $m_+(\lambda)=1$, the claim follows from the case (ii) of Lemma \ref{type D vanishing 300}.
For the general case, the induction hypothesis and the claim for the case (i) shows our claim by applying (\ref{type D vanishing 900}) to reduce the multiplicity for $\tau_{\A_{\rnk}}$ in $\tau_{\A_{1}}\cdots\tau_{\A_{\rnk}}$.
The claim for the case (iii) can be proved similarly.
\qed

\vspace{20pt}
Now, Proposition \ref{zigzag lemma type D} follows from Lemma \ref{type D vanishing 1200} and the case (iii) of Lemma \ref{type D vanishing 300}.

\vspace{10pt}

For a signed Young diagram $\lambda$ with $\rnk$ rows fitting into the $\rnk\times\rnk$-square, let
\begin{align*}
&m:=|\{i \mid \text{$\lambda_i=\rnk$} \}|=m_+(\lambda) + m_-(\lambda) 
\end{align*}
be the number of rows of $\lambda$ of length $\rnk$.
Recall that the numbers $\parasec_{r}, \parathr_{r}, \parafir_{r}$, and $\parafor_{r}$ are defined in (\ref{definition of a b c}) and (\ref{definition of y}).
We now define
\begin{align*}
\widetilde{\parafor}_{1}
:=
\begin{cases}
\displaystyle{
2^{(\rnk-\lambda_1-1)(1-m)}
\binom{ \parasec_{1} }{ \parafir_{1} } 
\binom{ \parathr_{1} }{ \parafir_{1} } 
}
\qquad &\text{if $m\leq1$}, \vspace{3pt} \\ 
\displaystyle{
-\binom{\parathr_1-1}{\parafir_1-1}
}
&\text{if $m\geq2$ and $m_+(\lambda)m_-(\lambda)\neq0$}, \vspace{3pt} \\
\displaystyle{
\left(2^{\parasec_1}-\parasec_1-1\right)
 \binom{ \parathr_{1} }{ \parafir_{1} }
 +
 \binom{ \parathr_{1}-1 }{ \parafir_{1} } 
}
&\text{if $m\geq2$ and $m_+(\lambda)m_-(\lambda)=0$}
\end{cases}
\end{align*}

\begin{theorem}\label{main thm for type D}
If $\emptyset\subsetneq \A_1,\cdots,\A_{\rnk}\subsetneq[\pm\rnk]$ satisfying \emph{(\ref{type D basics 100})} form a subchain of a chain in $\mathcal{C}$, then we have
\begin{align*} 
\intnum{\tau_{\A_1}\cdots\tau_{\A_{\rnk}}} = (-1)^{\rnk+s}\widetilde{y}_1y_2\cdots y_s.
\end{align*}
where $\mu_X$ is the fundamental homology class of $X$ and $\lambda$ is the signed Young diagram consisting of $|\A_1|,\cdots,|\A_{\rnk}|$ reordered as a weakly decreasing sequence and $m$ is the number of rows of $\lambda$ of length $\rnk$.
Otherwise, the intersection number is zero.
\end{theorem}

\begin{remark}
\emph{
In each case, the given number vanishes unless each step of the zigzag line of the lower-right corners of $\lambda$ crosses the dotted anti-diagonal. 
}
\end{remark}
\noindent
\textit{Proof of Theorem \ref{main thm for type D}.}
\\ \no
We compute the intersection number $J(\lambda):=(-1)^{\rnk-s}\intnum{\tau_{\A_1}\cdots\tau_{\A_{\rnk}}}$ with sign
where we can assume that each step of the zigzag line of the corners of $\lambda$ crosses the dotted anti-diagonal by Proposition \ref{zigzag lemma type D} and the remark above.
We first prove the case (i).
If $m_+(\lambda)=1$ and $m_-(\lambda)=0$ (or $m_+(\lambda)=1$ and $m_-(\lambda)=0$), then it follows that $J(\lambda) = \parafor_{1}\cdots \parafor_{s}$. This can be proved by induction similar to that used in the proof of Theorem \ref{intro formula of type A intersection} because of the separating properties (\ref{type D vanishing 140}).
So, let us consider the case of $m_+(\lambda)=m_-(\lambda)=0$.
In this case, Proposition \ref{zigzag lemma type D} shows that the separation rule (\ref{type D vanishing 130}) replaces the square $\tau_{\A_{\rnk}}^2$ in $\tau_{\A_1}\cdots\tau_{\A_{\rnk}}$ to
\begin{align*}
\sum_{\substack{k\in \A', \ \pm l\notin \A' \\ \A'\subsetneq \A_{\rnk}, \ |B'|\neq\rnk}} \tau_{\A'} \tau_{\A_{\rnk}}
+
\sum_{\substack{k\in \A', \ -l\in \A' \\ |\A'|=\rnk}} \tau_{\A'} \tau_{\A_{\rnk}}
\end{align*}
for some $k\in\A_{\rnk}$ and $\pm l\in\A_{\rnk}$
when we compute the intersection number $J(\lambda)$ with sign.
Namely, this replacement can be pictured as
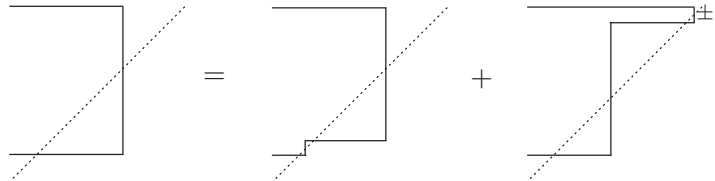
\[
\unitlength 0.1in
\begin{picture}( 36.2100, 10.1200)(  4.0000,-20.8800)
%
\special{pn 8}%
\special{pa 400 1170}%
\special{pa 992 1170}%
\special{fp}%
%
\special{pn 8}%
\special{pa 992 1946}%
\special{pa 400 1946}%
\special{fp}%
%
\special{pn 8}%
\special{pa 2368 1178}%
\special{pa 2368 1874}%
\special{fp}%
%
\special{pn 8}%
\special{pa 2368 1874}%
\special{pa 1950 1874}%
\special{fp}%
%
\special{pn 8}%
\special{pa 1946 1872}%
\special{pa 1946 1948}%
\special{fp}%
%
\special{pn 8}%
\special{pa 1774 1950}%
\special{pa 1950 1950}%
\special{fp}%
%
\special{pn 8}%
\special{pa 1774 1178}%
\special{pa 2368 1178}%
\special{fp}%
%
\special{pn 8}%
\special{pa 1314 1174}%
\special{pa 402 2086}%
\special{dt 0.030}%
%
\special{pn 8}%
\special{pa 3546 1256}%
\special{pa 3546 1952}%
\special{fp}%
%
\special{pn 8}%
\special{pa 3110 1950}%
\special{pa 3546 1950}%
\special{fp}%
%
\special{pn 8}%
\special{pa 3110 1174}%
\special{pa 3982 1174}%
\special{fp}%
%
\special{pn 8}%
\special{pa 3982 1256}%
\special{pa 3542 1256}%
\special{fp}%
%
\special{pn 8}%
\special{pa 4022 1174}%
\special{pa 3110 2086}%
\special{dt 0.030}%
\put(14.1300,-15.7000){\makebox(0,0)[lb]{$=$}}%
\put(28.1200,-15.9900){\makebox(0,0)[lb]{$+$}}%
%
\special{pn 8}%
\special{pa 992 1174}%
\special{pa 992 1946}%
\special{fp}%
%
\special{pn 8}%
\special{pa 2688 1178}%
\special{pa 1776 2088}%
\special{dt 0.030}%
%
\special{pn 8}%
\special{pa 3982 1256}%
\special{pa 3982 1178}%
\special{fp}%
\put(39.9600,-12.4600){\makebox(0,0)[lb]{$_{\pm}$}}%
\end{picture}%
\]
where we omit the coefficients in the picture.
Hence, with the claim for the previous case, we get $J(\lambda)=2^{\rnk-\lambda_1-1}\parafor_1\cdots \parafor_s$
as in the case of type $B_{\rnk}$.

Let us consider the case (ii-a). 
We prove the claim by induction on the sum of the multiplicities for $\A_i$'s satisfying $|\A_i|\neq\rnk$.
The base case has 
$\lambda_i=\rnk+1-i$ for all $\lambda_i\neq\rnk$, so it is obvious that $-J(\lambda)=(-1)^{\rnk-s-1}\intnum{\tau_{\A_1}\cdots\tau_{\A_{\rnk}}}$ is equal to $1$ by iterating (\ref{type D vanishing 1100}). For the general case, we apply (\ref{type D vanishing 140}) and (\ref{type D vanishing 141}) 
to some square ${\tau_{\A_i}}^2$ with $|\A_i|\neq\rnk$ in $J(\lambda)$. 
If $\coi_{2}<i$, the computation with (\ref{type D vanishing 140}) works as in the proof of Theorem \ref{intro formula of type A intersection}.
If $\coi_{1}<i \leq \coi_{2}$, then the right-hand-side of (\ref{type D vanishing 141}) applied to $-J(\lambda)$ can be calculated as follows by the induction hypothesis:
\begin{align*}
\binom{\parathr_2}{\parafir_2}
\cdot
\binom{\parathr_1-1}{\parafir_1-1}
\binom{\parasec_2-1}{\parafir_2-1}
\parafor_{3}\cdots \parafor_{s} 
+
\binom{\parathr_1-1}{\parafir_1-1}
\cdot
\binom{\parasec_2-1}{\parafir_2}\binom{\parathr_2}{\parafir_2}
\parafor_{3}\cdots \parafor_{s} 
\end{align*}
which is equal to $\binom{\parathr_1-1}{\parafir_1-1}
\parafor_{2}\cdots \parafor_{s}$, and the claim follows.

Let us consider the case (ii-b). We can assume $m=m_+(\lambda)(=\parasec_1+1)$ without loss of generality.
We first consider the case that $\lambda_i=\rnk+1-i$ for all $\lambda_i\neq\rnk$. Note that 
$\binom{\parathr_1-1}{\parafir_1}=0$
in this case.
For this special case, we prove the claim by induction on $m=m_+(\lambda)$.
For the case $m=\parasec_1+1=2$, the intersection number is zero by (\ref{type D vanishing 1000}), and the claim follows since $2^{\parasec_1}-\parasec_1-1=0$.
For the general case, we have an inductive formula by (\ref{type D vanishing 900}), namely, 
\[
\input{pic-typeD-3} 
\]
where the first coefficient $\parasec_1-1$ comes from the choices of an element in $\Dsetsec'\backslash(\Dsetfir\cup\{k,l\})$ turned to be negative, and the second coefficient 2 comes from the two summands corresponding to $k\in \Dsetsec', \pm l\notin \Dsetsec'$ and $\pm k\notin \Dsetsec', l\in \Dsetsec'$.
Noticing that the intersection number for the first summand is equal to $1$ as we already saw, it follows that 
\begin{align*}
J(\lambda)
= \sum_{i=0}^{\parasec_1} 2^i (\parasec_1-1-i) 
= 2^{\parasec_1} - \parasec_1 -1.
\end{align*}

We now prove the claim (ii-b) by induction on the sum of the multiplicities for $\A_i$ satisfying $|\A_i|\neq\rnk$.
The base case $\lambda_i=\rnk+1-i$ for all $\lambda_i\neq\rnk$ is proved above.
For the general case, we apply (\ref{type D vanishing 140}). 
If $\coi_{2}<i$, the computation with (\ref{type D vanishing 140}) again works as in the proof of Theorem \ref{intro formula of type A intersection}.
If $\coi_{1}<i \leq \coi_{2}$, we also apply (\ref{type D vanishing 140}) to a square ${\A_i}^2$ in $-J(\lambda)$
Namely, we have
\[
\unitlength 0.1in
\begin{picture}( 49.7900, 10.1400)(  8.0000,-20.8400)
%
\special{pn 8}%
\special{pa 800 1168}%
\special{pa 1678 1168}%
\special{fp}%
%
\special{pn 8}%
\special{pa 1678 1168}%
\special{pa 1678 1398}%
\special{fp}%
%
\special{pn 8}%
\special{pa 1678 1398}%
\special{pa 1240 1398}%
\special{fp}%
%
\special{pn 8}%
\special{pa 1240 1398}%
\special{pa 1240 1948}%
\special{fp}%
%
\special{pn 8}%
\special{pa 1240 1948}%
\special{pa 800 1948}%
\special{fp}%
\put(17.0000,-12.4000){\makebox(0,0)[lb]{$_+$}}%
\put(17.0000,-13.1700){\makebox(0,0)[lb]{$_+$}}%
\put(17.0000,-13.9400){\makebox(0,0)[lb]{$_+$}}%
%
\special{pn 8}%
\special{pa 2558 1398}%
\special{pa 2558 1870}%
\special{fp}%
%
\special{pn 8}%
\special{pa 2558 1870}%
\special{pa 2294 1870}%
\special{fp}%
%
\special{pn 8}%
\special{pa 2294 1870}%
\special{pa 2294 1948}%
\special{fp}%
%
\special{pn 8}%
\special{pa 2118 1948}%
\special{pa 2294 1948}%
\special{fp}%
%
\special{pn 8}%
\special{pa 2118 1168}%
\special{pa 2996 1168}%
\special{fp}%
%
\special{pn 8}%
\special{pa 2996 1168}%
\special{pa 2996 1398}%
\special{fp}%
%
\special{pn 8}%
\special{pa 2996 1398}%
\special{pa 2558 1398}%
\special{fp}%
\put(30.1000,-12.4000){\makebox(0,0)[lb]{$_+$}}%
\put(30.1000,-13.1700){\makebox(0,0)[lb]{$_+$}}%
\put(30.1000,-13.9400){\makebox(0,0)[lb]{$_+$}}%
%
\special{pn 8}%
\special{pa 3034 1168}%
\special{pa 2118 2084}%
\special{dt 0.030}%
%
\special{pn 8}%
\special{pa 3874 1476}%
\special{pa 3874 1948}%
\special{fp}%
%
\special{pn 8}%
\special{pa 3436 1948}%
\special{pa 3874 1948}%
\special{fp}%
%
\special{pn 8}%
\special{pa 3436 1168}%
\special{pa 4314 1168}%
\special{fp}%
%
\special{pn 8}%
\special{pa 4314 1168}%
\special{pa 4314 1398}%
\special{fp}%
%
\special{pn 8}%
\special{pa 4314 1398}%
\special{pa 4084 1398}%
\special{fp}%
\put(43.3000,-12.4000){\makebox(0,0)[lb]{$_+$}}%
\put(43.3000,-13.1700){\makebox(0,0)[lb]{$_+$}}%
\put(43.3000,-13.9400){\makebox(0,0)[lb]{$_+$}}%
%
\special{pn 8}%
\special{pa 3874 1476}%
\special{pa 4084 1476}%
\special{fp}%
%
\special{pn 8}%
\special{pa 4084 1476}%
\special{pa 4084 1398}%
\special{fp}%
%
\special{pn 8}%
\special{pa 1718 1168}%
\special{pa 800 2084}%
\special{dt 0.030}%
%
\special{pn 8}%
\special{pa 4352 1168}%
\special{pa 3436 2084}%
\special{dt 0.030}%
%
\special{pn 8}%
\special{pa 4862 1168}%
\special{pa 5742 1168}%
\special{fp}%
%
\special{pn 8}%
\special{pa 5742 1168}%
\special{pa 5742 1476}%
\special{fp}%
%
\special{pn 8}%
\special{pa 5742 1476}%
\special{pa 5302 1476}%
\special{fp}%
%
\special{pn 8}%
\special{pa 5302 1476}%
\special{pa 5302 1948}%
\special{fp}%
%
\special{pn 8}%
\special{pa 5302 1948}%
\special{pa 4862 1948}%
\special{fp}%
\put(57.6000,-12.4000){\makebox(0,0)[lb]{$_+$}}%
\put(57.6000,-13.1700){\makebox(0,0)[lb]{$_+$}}%
\put(57.6000,-13.9400){\makebox(0,0)[lb]{$_+$}}%
%
\special{pn 8}%
\special{pa 5780 1168}%
\special{pa 4862 2084}%
\special{dt 0.030}%
\put(57.6000,-14.7100){\makebox(0,0)[lb]{$_-$}}%
\put(18.2000,-16.5000){\makebox(0,0)[lb]{$=$}}%
\put(31.3000,-16.5000){\makebox(0,0)[lb]{$+$}}%
\put(45.6000,-16.5000){\makebox(0,0)[lb]{$+$}}%
\end{picture}%
\]
with omitting the coefficients.
The right-hand-side can be calculated by the induction hypothesis and the claim for the case (ii-a), and the intersection number $J(\lambda)=(-1)^{\rnk-s}\intnum{\tau_{\A_1}\cdots\tau_{\A_{\rnk}}}$ is \begin{align*}
&
\binom{\parathr_2}{\parafir_2} \cdot
\Bigg\{
(2^{\parasec_1}-\parasec_1-1)
\binom{\parathr_1}{\parafir_1}
\binom{\parasec_2-1}{\parafir_2-1}
\parafor_{3}\cdots \parafor_{s} 
+
\binom{\parathr_1-1}{\parafir_1}
\binom{\parasec_2-1}{\parafir_2-1}
\parafor_3\cdots \parafor_s \Bigg\}
\\
&\qquad\quad
+ 
\binom{\parathr_1}{\parafir_1} \cdot
\Bigg\{
(2^{\parasec_1}-\parasec_1-1)
\binom{\parasec_2-1}{\parafir_2}\binom{\parathr_2}{\parafir_2}
\parafor_{3}\cdots \parafor_{s} 
+ \ 
0 \ 
\Bigg\}
\\
&\qquad\quad
+
\binom{\parathr_1-1}{\parafir_1+1-1}
\binom{\parasec_2-1}{\parafir_2}\binom{\parathr_2}{\parafir_2}
\parafor_3\cdots \parafor_s \\
&\quad
=
(2^{\parasec_1}-\parasec_1-1)
\binom{\parathr_1}{\parafir_1}
\parafor_{2}
\cdots \parafor_{s} 
+
\binom{\parathr_1-1}{\parafir_1}
\parafor_2\cdots \parafor_s.
\end{align*}
\qed
\\ 

For example, for $\rnk=5$ with the convention $\bar{k}=-k$, we can calculate
\begin{align*}
\intnum{{ \tau_{\{\bar{1}\}} }^2 { \tau_{\{\bar{1},3,4,5,\bar{2}\}} }^3}  = -4
\end{align*}
by the case (ii-b) of Theorem \ref{main thm for type D}. See Figure \ref{EX for type D} for the Young diagram corresponding to ${ \tau_{\{\bar{1}\}} }^2 { \tau_{\{\bar{1},3,4,5,\bar{2}\}} }^3$.
\begin{figure}[h]
\centering
\unitlength 0.1in
\begin{picture}(  7.3000,  8.9000)( 24.6000,-19.8000)
\put(26.3000,-12.6000){\makebox(0,0)[lb]{$3$}}%
\put(27.7000,-12.6000){\makebox(0,0)[lb]{$4$}}%
\put(29.1000,-12.6000){\makebox(0,0)[lb]{$5$}}%
\put(30.5000,-12.6000){\makebox(0,0)[lb]{$\bar{2}$}}%
%
\special{pn 8}%
\special{pa 3160 1420}%
\special{pa 3020 1420}%
\special{pa 3020 1280}%
\special{pa 3160 1280}%
\special{pa 3160 1420}%
\special{fp}%
%
\special{pn 8}%
\special{pa 3160 1560}%
\special{pa 3020 1560}%
\special{pa 3020 1420}%
\special{pa 3160 1420}%
\special{pa 3160 1560}%
\special{fp}%
%
\special{pn 8}%
\special{pa 3020 1560}%
\special{pa 2880 1560}%
\special{pa 2880 1420}%
\special{pa 3020 1420}%
\special{pa 3020 1560}%
\special{fp}%
%
\special{pn 8}%
\special{pa 3020 1420}%
\special{pa 2880 1420}%
\special{pa 2880 1280}%
\special{pa 3020 1280}%
\special{pa 3020 1420}%
\special{fp}%
%
\special{pn 8}%
\special{pa 2880 1420}%
\special{pa 2740 1420}%
\special{pa 2740 1280}%
\special{pa 2880 1280}%
\special{pa 2880 1420}%
\special{fp}%
%
\special{pn 8}%
\special{pa 2880 1560}%
\special{pa 2740 1560}%
\special{pa 2740 1420}%
\special{pa 2880 1420}%
\special{pa 2880 1560}%
\special{fp}%
%
\special{pn 8}%
\special{pa 2740 1560}%
\special{pa 2600 1560}%
\special{pa 2600 1420}%
\special{pa 2740 1420}%
\special{pa 2740 1560}%
\special{fp}%
%
\special{pn 8}%
\special{pa 2740 1420}%
\special{pa 2600 1420}%
\special{pa 2600 1280}%
\special{pa 2740 1280}%
\special{pa 2740 1420}%
\special{fp}%
%
\special{pn 8}%
\special{pa 2600 1560}%
\special{pa 2460 1560}%
\special{pa 2460 1420}%
\special{pa 2600 1420}%
\special{pa 2600 1560}%
\special{fp}%
%
\special{pn 8}%
\special{pa 2600 1420}%
\special{pa 2460 1420}%
\special{pa 2460 1280}%
\special{pa 2600 1280}%
\special{pa 2600 1420}%
\special{fp}%
%
\special{pn 8}%
\special{pa 2600 1700}%
\special{pa 2460 1700}%
\special{pa 2460 1560}%
\special{pa 2600 1560}%
\special{pa 2600 1700}%
\special{fp}%
%
\special{pn 8}%
\special{pa 2600 1840}%
\special{pa 2460 1840}%
\special{pa 2460 1700}%
\special{pa 2600 1700}%
\special{pa 2600 1840}%
\special{fp}%
%
\special{pn 8}%
\special{sh 0.300}%
\special{pa 2600 1980}%
\special{pa 2460 1980}%
\special{pa 2460 1840}%
\special{pa 2600 1840}%
\special{pa 2600 1980}%
\special{fp}%
%
\special{pn 8}%
\special{pa 2740 1700}%
\special{pa 2600 1700}%
\special{pa 2600 1560}%
\special{pa 2740 1560}%
\special{pa 2740 1700}%
\special{fp}%
%
\special{pn 8}%
\special{pa 2880 1700}%
\special{pa 2740 1700}%
\special{pa 2740 1560}%
\special{pa 2880 1560}%
\special{pa 2880 1700}%
\special{fp}%
%
\special{pn 8}%
\special{pa 3020 1700}%
\special{pa 2880 1700}%
\special{pa 2880 1560}%
\special{pa 3020 1560}%
\special{pa 3020 1700}%
\special{fp}%
%
\special{pn 8}%
\special{sh 0.300}%
\special{pa 3160 1700}%
\special{pa 3020 1700}%
\special{pa 3020 1560}%
\special{pa 3160 1560}%
\special{pa 3160 1700}%
\special{fp}%
%
\special{pn 8}%
\special{pa 2460 1980}%
\special{pa 3160 1280}%
\special{fp}%
\put(24.9000,-12.6000){\makebox(0,0)[lb]{$\bar{1}$}}%
\put(31.9000,-14.0000){\makebox(0,0)[lb]{$+$}}%
\put(31.9000,-15.4000){\makebox(0,0)[lb]{$+$}}%
\put(31.9000,-16.8000){\makebox(0,0)[lb]{$+$}}%
\end{picture}%
\caption{}
\label{EX for type D}
\end{figure}

\vspace{10pt}
For each even signed permutation $u\in \widetilde{\mathfrak{S}}_{\rnk}^+$, 
an element $i\in[\rnk]$ satisfies $u(\alpha_i)\in\Phi^-$ if and only if
\begin{itemize}
 \item[(D-1)] if $i\leq\rnk-2$, then $u(i)>u(i+1)$ with the same sign, or $u(i)<u(i+1)$ with different signs,
 \item[(D-2)] if $i=\rnk-1$, then $u(\rnk-1), u(\rnk)<0$, or 
$u(\rnk-1)$ and $u(\rnk)$ have different signs and the absolute value of the negative one is less than the positive one,
 \item[(D-3)] if $i=\rnk$, then $u(\rnk-1)>u(\rnk)$ with the same sign, or $u(\rnk-1)<u(\rnk)$ with different signs.
\end{itemize}
Consider the similar condition 
\begin{itemize}
 \item[(A-1)] if $i\leq n-2$, then $u(i)<u(i+1)$ with the same sign or $u(i)>u(i+1)$ with different signs,
 \item[(A-2)] if $i=\rnk-1$, then $u(\rnk-1), u(\rnk)>0$, or 
$u(\rnk-1)$ and $u(\rnk)$ have different signs and the absolute value of the negative one is greater than the positive one,
 \item[(A-3)] if $i\neq n$, then $u(\rnk-1)<u(\rnk)$ with the same sign or $u(\rnk-1)>u(\rnk)$ with different signs,
\end{itemize}
Denote
\begin{align*}
D(u):=&\{ u[i] \mid \text{$i\leq\rnk-2$ and $i$ satisfies (D)} \} \\
&\qquad
\cup \{ u[\rnk]_+ \mid \text{$i=\rnk-1$ satisfies (D)}\}
\cup \{ u[\rnk]_- \mid \text{$i=\rnk$ satisfies (D)}\} \\
A(u):=&\{ u[i] \mid \text{$i$ satisfies (A)} \} \\
&\qquad
\cup \{ u[\rnk]_+ \mid \text{$i=\rnk-1$ satisfies (A)}\}
\cup \{ u[\rnk]_- \mid \text{$i=\rnk$ satisfies (A)}\}
\end{align*}
where 
\begin{align*}
[\rnk]_+ = \{1,2,\cdots,\rnk-1,\rnk\} \quad \text{and} \quad
[\rnk]_- = \{1,2,\cdots,\rnk-1,-\rnk\}
\end{align*}
(cf. Figure \ref{type D chain}).
We define a signed Young diagram $\lambda_{u,v}^w$ for $u,v,w\in \widetilde{\mathfrak{S}}_{\rnk}^+$ in the manner described in  Section \ref{ring str of coh}.
Note that we put $I(\emptyset)=0$ as a convention.

Now, the intersection number of $Y^w$, $X_u$ and $X_v$ in $X$ of type $D_{\rnk}$ is given by the following. 
\begin{corollary}\label{triple intersection for type D}
For even signed permutations $u,v,w\in \widetilde{\mathfrak{S}}_{\rnk}^+$, we have
\[ \intnum{[Y^w][X_u][X_v]}=\widetilde{I}(\WY{u}{v}{w}) \]
where $\widetilde{I}=(-1)^{\rnk+s}\widetilde{y}_1y_2\cdots y_s$ is the function described in \emph{Theorem \ref{main thm for type D}}.
\end{corollary}

\vspace{10pt}
For example, for $\rnk=5$ with the convention $\bar{k}=-k$, the Young diagram $\lambda_{\bar{1}345\bar{2}, \bar{1}345\bar{2}}^{\bar{1}\bar{2}543}$ is the one in Figure \ref{EX for type D}, and hence we obtain
\begin{align*}
\intnum{[Y^{\bar{1}\bar{2}543}][X_{\bar{1}345\bar{2}}][X_{\bar{1}345\bar{2}}]}  = -4.
\end{align*}

\section{On exceptional types}
In this section, we include the computation of intersection numbers of invariant divisors of the toric manifold $X$ for the root system of exceptional type $G_2$.
For other exceptional types $F_4, E_6, E_7$, and $E_8$, 
it would be interesting to find combinatorial objects which effectively compute the intersection numbers of invariant divisors.\\


Let $E=\{x\in\R^{3} \mid x_1+x_2+x_3=0\}$.
The roots are 
\begin{align*}
 &\pm(t_1-t_2), \ \pm(t_1-t_3), \ \pm(t_2-t_3), \\
 &\pm(2t_1-t_2-t_3), \ \pm(2t_2-t_1-t_3), \ \pm(2t_3-t_1-t_2).
\end{align*} 
where $t_i\in\R^{3}$ is the $i$-th standard vector.
We choose $\SimR=\{t_1-t_2, -2t_1+t_2+t_3\}$
as the set of simple roots, and write $\alpha_1=t_1-t_2$ and $\alpha_2=-2t_1+t_2+t_3$.
The Weyl group $W$ is the dihedral group of order 12 which is identified with the subgroup
\[
\WG:=
\{ u\in \widetilde{\mathfrak{S}}_3 \mid \text{$u(1)$, $u(2)$, and $u(3)$ have the same sign} \}
\]
of the 3rd signed permutation group.
Under this identification, the action of the Weyl group on $E$ is written as the natural action of $\WG$ on the indexes $i$ of $t_i$; $u\cdot t=t_{u(1)}$ for $u\in \WG$ where $t_{-i}:=-t_i$ ($1\leq i\leq 3$). 
This action of $W_{G_2}$ preserves $\Phi$.
The minimal generators $\fcw_1, \fcw_2\in E^*$ of the fundamental Weyl chamber $\sigma_{\text{id}}$ are
\[
\fcw_1=e_3-e_2, \ \fcw_2=\frac{1}{3}(2e_3-e_1-e_2)
\]
where $\{e_i\}_i\subset (\R^3)^*$ is the dual basis of $\{t_i\}_i\subset \R^3$.

Denoting by $2^{[\pm3]}$ the set of all subsets of $[\pm3]=\{1,2,3,-1,-2,-3\}$, we have a well-defined map $\dPhi \rightarrow 2^{[\rnk+1]}$ by sending 
\[
e_{u(3)}-e_{u(2)} \mapsto \{u(3), -u(2)\}, \quad
\frac{1}{3}(2e_{u(3)}-e_{u(1)}-e_{u(2)}) \mapsto \{u(3)\}
\]
for $u\in \WG$.
This is an injection, and hence we can identify $\Phi^*$ with the following subset of $2^{[\pm3]}$; 
\begin{align*}
\mathcal{S}:=\{ 3\bar{2}, \bar{3}2, 3\bar{1}, \bar{3}1, 2\bar{1}, \bar{2}1, 3,\bar{3}, 2, \bar{2}, 1, \bar{1}\}
\end{align*} 
where $\bar{k}=-k$ for $1\leq k\leq3$ and each sequence $ab$ in $\mathcal{S}$ is the set $\{a,b\}$, i.e. $3\bar{2}=\{3,-2\}$ for example.
Now, for each $S\in \mathcal{S}$, we have 
$\tau_S:=\tau_{u\omega_i}\in H^2(X)$ where $u\omega_i\in\dPhi$ corresponds to $S$ by this identification. 
Then, for $S_1,S_2\in\mathcal{S}$, it follows by Lemma \ref{prelim rays generating a cone} that $\tau_{S_1}\tau_{S_2}=0$ unless these sets form a nested chain of subsets, i.e. $S_1\subset S_2$ or $S_1\supset S_2$.

The linear relations (\ref{prelim linear relation}) for $\alpha=\alpha_1,\alpha_2$ are translated to 
\begin{align*}
&\tau_{3\bar{2}} 
+ \tau_{\bar{3}1} 
+ 2\tau_{\bar{2}1}  
+ \tau_{\bar{2}}  
+ \tau_{1}
=
\tau_{\bar{3}2}
+\tau_{3\bar{1}} 
+2\tau_{2\bar{1}} 
+\tau_{2}  
+ \tau_{\bar{1}}, \\
&3\tau_{3\bar{1}} 
+ 3\tau_{2\bar{1}}  
+ \tau_{3} 
+ \tau_{2} 
+ 2\tau_{\bar{1}} 
=
3\tau_{\bar{3}1}
+3\tau_{\bar{2}1} 
+\tau_{\bar{3}}
+\tau_{\bar{2}}   
+2\tau_{1},
\end{align*}
respectively.
From these relations together with the above observation about the vanishing of  $\tau_{S_1}\tau_{S_2}$, we see that
\begin{align*}
\tau_{3\bar{2}}\tau_{3}=1, \ 
\tau_{3\bar{2}}\tau_{3\bar{2}}=-1, \ 
\tau_{3}\tau_{3}=-3.
\end{align*}
Now, let 
\[
I_{G_2}(2,1):=1, \ I_{G_2}(1,1):=-3, \ I_{G_2}(2,2):=-1.
\]
where $(2,1)$, $(1,1)$, and $(2,2)$ are Young diagrams with 2 rows.
Now the next claim follows from Lemma \ref{prelim Weyl inv}; 
if $S_1,S_2\in\mathcal{S}$ form a nested chain of subsets, 
then we have
\begin{align}\label{type G2 int num}
\intnum{\tau_{\A_1}\tau_{\A_2}} 
= I_{G_2}(\lambda)
\end{align}
where $\mu_X$ is the fundamental homology class and $\lambda$ is the Young diagram consisting of $|\A_1|$ and $|\A_2|$ reordered as a weakly decreasing sequence. Otherwise, the intersection number is zero.

Finally, we list the presentations of $[X_u]$ as monomials of $\tau_{\A}$ for all $u\in \WG$ in one-line notations;
\begin{align*}
&
[X_{123}]=1, \hspace{21pt}
[X_{213}]=\tau_{3\bar{1}}, 
[X_{132}]=\tau_{2\bar{3}}, 
[X_{231}]=\tau_{1\bar{3}}, 
[X_{312}]=\tau_{2}, \hspace{7pt}
[X_{321}]=\tau_{1}, \\
&
[X_{\bar{1}\bar{2}\bar{3}}]=\tau_{\bar{3}2}\tau_{\bar{3}}, \ 
[X_{\bar{2}\bar{1}\bar{3}}]=\tau_{\bar{3}}, \ 
[X_{\bar{1}\bar{3}\bar{2}}]=\tau_{\bar{2}}, \ \hspace{1pt}
[X_{\bar{2}\bar{3}\bar{1}}]=\tau_{\bar{1}}, \ 
[X_{\bar{3}\bar{1}\bar{2}}]=\tau_{\bar{2}1}, \ 
[X_{\bar{3}\bar{2}\bar{1}}]=\tau_{\bar{1}2}.
\end{align*}
Since we have $[Y^u]=(w_0^{-1})^*[X_{w_0 u}]=w_0^*[X_{w_0 u}]$ where $w_0=\bar{1}\bar{2}\bar{3}$ is the longest permutation, we obtain the list of  $[Y^u]$;
\begin{align*}
&
[Y^{\bar{1}\bar{2}\bar{3}}]=1, \hspace{22pt}
[Y^{\bar{2}\bar{1}\bar{3}}]=\tau_{\bar{3}1}, 
[Y^{\bar{1}\bar{3}\bar{2}}]=\tau_{\bar{2}3}, 
[Y^{\bar{2}\bar{3}\bar{1}}]=\tau_{\bar{1}3}, 
[Y^{\bar{3}\bar{1}\bar{2}}]=\tau_{\bar{2}}, \hspace{7pt}
[Y^{\bar{3}\bar{2}\bar{1}}]=\tau_{\bar{1}}, \\
&
[Y^{123}]=\tau_{3\bar{2}}\tau_{3}, \ 
[Y^{213}]=\tau_{3}, \ \hspace{1pt}
[Y^{132}]=\tau_{2}, \ 
[Y^{231}]=\tau_{1}, \ 
[Y^{312}]=\tau_{2\bar{1}}, \ 
[Y^{321}]=\tau_{1\bar{2}}.
\end{align*}
With these lists, we can compute intersection numbers $\intnum{[Y^w][X_u][X_v]}$ for all $u,v,w\in \WG$ by (\ref{type G2 int num}).

\end{document}